\documentclass[12pt]{amsart}
\usepackage{a4wide,amsmath,enumerate,xcolor,graphicx}
\usepackage{comment,enumitem}
\allowdisplaybreaks

\let\pa\partial
\let\na\nabla
\let\eps\varepsilon
\newcommand{\N}{{\mathbb N}}
\newcommand{\R}{{\mathbb R}}
\newcommand{\diver}{\operatorname{div}}

\newcommand{\dom}{\mathcal{D}}
\newcommand{\sym}{\operatorname{sym}}
\newcommand{\D}{\mathrm{D}}
\newcommand{\dd}{\mathrm{d}}

\newtheorem{theorem}{Theorem}[section]
\newtheorem{lemma}[theorem]{Lemma}
\newtheorem{proposition}[theorem]{Proposition}
\theoremstyle{remark}
\newtheorem{remark}[theorem]{Remark}

\theoremstyle{plain}
\newtheorem{example}[theorem]{Example}

\begin{document}

\title[Construction of entropies]{Construction of entropies for cross-diffusion systems}

\author[X. Chen]{Xiuqing Chen}
\address{School of Mathematics (Zhuhai), Sun Yat-sen University, Zhuhai 519082, Guangdong Province, China}
\email{chenxiuqing@mail.sysu.edu.cn}

\author[A. J\"ungel]{Ansgar J\"ungel}
\address{Institute of Analysis and Scientific Computing, TU Wien, Wiedner Hauptstra\ss e 8--10, 1040 Wien, Austria}
\email{juengel@tuwien.ac.at} 

\author[M. Zhang]{Mingyue Zhang}
\address{Institute of Analysis and Scientific Computing, TU Wien, Wiedner Hauptstra\ss e 8--10, 1040 Wien, Austria}
\email{mingyue.zhang@tuwien.ac.at} 

\date{\today}

\thanks{The first author acknowledges support from the National
Natural Science Foundation of China (NSFC), grant 12471206. The last two authors acknowledge partial support from the Austrian Science Fund (FWF), grant 10.55776/PAT2687825, and from the Austrian Federal Ministry for Women, Science and Research and implemented by \"OAD, project MULT09/2025. This work has received funding from the European Research Council (ERC) under the European Union's Horizon 2020 research and innovation programme, ERC Advanced Grant NEUROMORPH, no.~101018153. For open-access purposes, the authors have applied a CC BY public copyright license to any author-accepted manuscript version arising from this submission.} 

\begin{abstract}
The analysis of cross-diffusion systems with nonsymmetric and indefinite diffusion matrices often relies on the existence of entropies, that is, Lyapunov functionals whose dissipation often yields gradient estimates. We propose a general strategy for constructing entropies for two-species cross-diffusion systems, based on the assumption that the mobility matrix is symmetric, in accordance with Onsager's reciprocal relations. The symmetry condition leads to a linear hyperbolic equation in two variables, which can be reformulated as a linear transport system along the eigendirections of the diffusion matrix, related to the approach of Jenssen and Kogan for hyperbolic conservation laws. Positive solutions to this system yield a family of entropies parametrized by two arbitrary functions of one variable. Sufficient conditions for their existence on (possibly unbounded) characteristic rectangles are given, the mechanisms that may prevent the existence of an entropy are identified, and systems without any entropy structure are presented. The general construction recovers known entropies for models arising in physics and biology and provides new ones. For systems with an arbitrary number of components, we derive sufficient conditions for the existence of separable entropies and give a rank criterion for the uniqueness of entropies (up to affine functions). Cross-diffusion systems without symmetry of the mobility matrix are also discussed.
\end{abstract}

% \paragraph{Keywords:}  
% \paragraph{AMS classification:}  

\maketitle

%%%%%%%%%%%%%%%%%%%%%%%%%%%%%%%%%%%%%%%%%%%%%%%%%%%%%%%%%%%%%%%%%%%

\section{Introduction}

Cross-diffusion systems arise naturally in the modeling of diffusive transport in multicomponent systems, with applications ranging from population dynamics and fluid mixtures to biological processes at the cellular level. Characteristic features of such systems are that the flux of each component depends on the density gradients of all interacting components and that the diffusion matrix is generally neither symmetric nor positive definite. This substantially complicates the mathematical analysis and, in particular, the derivation of a priori estimates for its solutions.

Motivated by thermodynamic principles, entropy functionals often provide a powerful tool for obtaining such estimates. However, identifying an appropriate entropy can be highly delicate, especially for diffusion models whose thermodynamic structure is not explicit or that arise beyond the classical framework of nonequilibrium thermodynamics. In this paper, we establish necessary and sufficient conditions for the existence of entropy functionals and provide a constructive procedure for determining them in two-species systems (if they exist). In particular, we show that a cross-diffusion system may admit no entropy, a unique entropy up to affine functions, or infinitely many distinct entropies. For some applications, we explain how the thermodynamically motivated entropy or free energy can be identified among the infinitely many admissible entropies.

We consider cross-diffusion systems of the form
\begin{align}\label{1.eq}
  \pa_t u_i = \diver\bigg(\sum_{j=1}^n A_{ij}(u)\na u_j\bigg)
  \quad\mbox{in }\Omega,\ t>0,\ i=1,\ldots,n,
\end{align}
supplemented with initial and no-flux boundary conditions, where $u_i=u_i(x,t)$ is the density of the $i$th component, $u=(u_1,\ldots,u_n)$ is the solution vector, $A_{ij}(u)$ are the diffusion coefficients, and $\Omega\subset\R^d$ is a bounded domain. The functions are assumed to have values in the domain $\dom\subset\R^n$. Examples are $\dom=(0,\infty)^n$ ($u_i$ is a density) and $\dom=\{u\in(0,1)^n:\sum_{i=1}^nu_i<1\}$ ($u_i$ is a volume fraction).

A minimal condition for local solvability is the positive stability of $A(u)$, meaning that all eigenvalues of $A(u)$ have positive real parts \cite[Sec.~4]{Ama93}. Establishing global-in-time existence, however, requires suitable a priori estimates. One powerful approach to obtaining such estimates is provided by the entropy method \cite{Jue15}. For this, let $h:\dom\to\R$ be a smooth convex function with Hessian matrix $\na^2 h$. A formal computation shows that
\begin{align*}
  \frac{\dd}{\dd t}\int_\Omega h(u)\dd x 
  + \sum_{i,j=1}^n\int_\Omega (\na^2hA)_{ij}(u)
  \na u_i\cdot\na u_j\dd x = 0.
\end{align*}
Then $t\mapsto\int_\Omega h(u(t))\dd x$ is a Lyapunov functional if the symmetric part of $\na^2hA$ is positive semidefinite. Moreover, if the symmetric part is positive definite, we may conclude a priori estimates for $\na u_i$ \cite{Jue15}. Thus, our task is to construct a convex function $h$ such that the symmetric part of $\na^2 h(u)A(u)$ is positive definite for all $u\in\dom$. Such a function is called an {\em entropy density}.

The construction is related to nonequilibrium thermodynamics. Indeed, introducing the so-called entropy variables (or chemical potentials) $w_i=\pa_ih:=\pa h/\pa u_i$, we have
\begin{align*}
  \frac{\dd}{\dd t}\int_\Omega h(u)\dd x 
  + \sum_{i,j=1}^n\int_\Omega B_{ij}(u)
  \na w_i\cdot\na w_j \dd x = 0,
\end{align*}
where $B(u) = A(u)(\na^2 h(u))^{-1}$ is referred to as the mobility (or Onsager) matrix, which is symmetric and positive definite in thermodynamic applications, in accordance with Onsager's reciprocal relations. In terms of the entropy variable $w$, we can formulate \eqref{1.eq} as
\begin{align*}
  \pa_t u_i - \sum_{j=1}^n\diver(B_{ij}(u)\na w_j) = 0
  \quad\mbox{in }\Omega,\ t>0,\ i=1,\ldots,n,
\end{align*}
where $u$ is interpreted as a function of $w$, $u(w)=(h')^{-1}(w)$ (assuming that $(h')^{-1}$ exists on $\R^n$).

The key idea of the paper is to exploit the symmetry requirement, i.e., we wish to construct entropies such that the mobility matrix is symmetric. (We also discuss cases in which the symmetry does not hold; see Section \ref{sec.sg}.) This requirement is related to the symmetrization of hyperbolic conservation laws, and therefore we call an entropy density $h$ a {\em convex symmetrizer} if $\na^2hA$ is symmetric. 

For $2\times 2$ matrices, this symmetry is formulated as the single equation
\begin{align}\label{1.hyper}
  A_{12}(u)\pa_{11}^2 h + (A_{22}(u)-A_{11}(u))\pa_{12}^2 h
  - A_{21}(u)\pa_{22}^2 h = 0 \quad\mbox{in }\dom,
\end{align}
where $\pa_{ij}^2=\pa^2/\pa u_i\pa u_j$. We show in Section \ref{sec.prep} that the strict convexity of $h$ and the symmetry of $\na^2hA$ imply that \eqref{1.hyper} is strictly hyperbolic (except in a trivial case). We reduce this equation in Section \ref{sec.main} to a first-order transport system for the diagonal entries of the Hessian $\na^2h$, formulated in the eigendirections of $A$. 

More precisely, let $r_\pm=r_\pm(u)$ be the (linearly independent) eigenvectors of $A(u)$. The symmetry of $\na^2hA$ implies that $r_\pm^T(\na^2h)r_\mp=0$, and we can construct the entropy density from the numbers $p_\pm:=r_\pm^T(\na^2 h)r_\pm$. The characteristic curves of \eqref{1.hyper} are the integral curves of the eigenvectors $r_\pm$. This turns \eqref{1.hyper} into a first-order system in characteristic form: $p_+$ is differentiated only along $r_-$, and $p_-$ only along $r_+$, leading to a so-called Goursat problem  \cite[Chap.~6]{RuRu94}. The Riemann invariants $\alpha_\pm$ straighten the characteristics, reducing the transport system to the normal form. Integrating factors remove the diagonal coefficients and leave the form $\pa_{\alpha_-}P_+ = K_+P_-$ and $\pa_{\alpha_+}P_- = K_-P_+$, where $\pa_{\alpha_\pm}$ are derivatives along $\alpha_\pm$ and $K_\pm$ are some coefficients. Then there is a convex symmetrizer if and only if {\em positive} solutions $P_\pm$ or $p_\pm$ exist, since a quadratic form without cross terms is positive definite if and only if its diagonal entries are positive; here, these entries are $p_\pm$.

This strategy, if applicable, gives two free functions of one variable. Hence, we obtain a family of entropy densities parametrized by two arbitrary functions of one variable. However, the method may fail. Possible failure mechanisms are the degeneracy of the eigenbasis (the eigenvectors become linearly dependent) and the failure of the characteristic curves to generate a global Goursat rectangle. Examples illustrating such situations are discussed in Section \ref{sec.fail}. Accordingly, several cases may arise: There may exist infinitely many convex symmetrizers, a unique convex symmetrizer up to affine functions, or no convex symmetrizer at all. 

We have made two essential assumptions: two variables $n=2$ and symmetry of the mobility matrix. The above strategy cannot generally be applied for more than two variables ($n>2$). The reason is that the problem is overdetermined; see Section \ref{sec.n} for details. In the case $n>2$, we exploit the total symmetry of the Hessian $\na^2h(u)$, yielding first-order equations for some coordinates, which split into a differential part (transport equations along the eigendirections) and an algebraic part, the latter being absent when $n=2$. We are able to provide in case $n>2$ sufficient conditions for the existence of {\em separable} convex symmetrizers, i.e., symmetrizers of the form $h(u)=\sum_{i=1}^n h_i(u_i)$, a condition on the rank of the algebraic part to obtain uniqueness of convex symmetrizers, as well as a criterion that rules out the existence of any entropy structure. When we relax the symmetry condition, we need to analyze matrix inequalities, whose general solution is an open problem. Some necessary conditions have been given in \cite{ChJu21}, based on matrix factorization theory. 
Our main results can be sketched as follows:

\begin{itemize}
\item $n=2$: A convex symmetrizer exists if and only if a linear transport system along the eigendirections of $A$ has positive solutions (necessary and sufficient condition); see Theorem \ref{thm.first}. This yields a family of convex symmetrizers parametrized by at most two free functions of one variable.
\item $n=2$: If the curvatures of the characteristic curves each have a constant sign and the domain $\dom$ is covered by one characteristic rectangle, then a convex symmetrizer exists (sufficient condition); see Theorem \ref{thm.const}.
\item $n=2$: If the coefficients of the transport system in normal form are sufficiently small, then a convex symmetrizer exists even when the curvatures change sign (sufficient condition); see Theorem \ref{thm.small}. 
\item $n=2$: If both characteristic curves are straight, a global convex symmetrizer can be constructed on a convex state domain without a characteristic-rectangle assumption (sufficient condition); see Theorem \ref{thm.straight}.
\item $n=2$: There exist matrices that do not have an entropy structure; see the example in \eqref{wind}, Proposition \ref{prop.nonex}, and Remark \ref{rem.nonex}. A counterexample in case $n=3$ is presented in Example \ref{ex.nonex}.
\item $n\ge 2$: We give an exact classification of the matrices $A$ that admit {\em separable} convex symmetrizers (necessary and sufficient conditions); see Theorem \ref{thm.class} for $n=2$. A sufficient condition for $n>2$ is proved in Theorem \ref{thm.classn}.
\item $n\ge 2$: When the quotients $A_{ij}/A_{ji}$ (for $A_{ji}\neq 0$) satisfy some easily checkable conditions, there exists a separable convex symmetrizer; see Theorem \ref{thm.septest} for $n=2$ and Theorem \ref{thm.septestn} for $n>2$. 
\item $n\ge 2$: The rank of the algebraic part determines the maximal number of independent convex symmetrizers. When the rank equals $n-1$, the convex symmetrizer is unique, up to a positive factor and affine functions (Theorem \ref{thm.unique}).
\end{itemize}

Generally, there exist families of entropies parametrized by at most $n$ free functions of one variable, or there is a unique convex symmetrizer (up to a factor and affine functions), or there exists no convex symmetrizer at all. We apply the results to some cross-diffusion systems arising in physics and biology:
\begin{itemize}
\item Shigesada--Kawasaki--Teramoto (SKT) population models (Examples \ref{ex.skt2} and \ref{ex.sktn}),
\item Busenberg--Travis-type models for populations of several species (Section \ref{sec.bt2}, Example \ref{ex.btn}),
\item two-species Maxwell--Stefan models for gas mixtures (Section \ref{sec.ms2}),
\item segregating granular material model (Section \ref{sec.sg}),
\item thin-film solar-cell model (Section \ref{sec.tfsc}).
\end{itemize}
In particular, we construct new convex symmetrizers for the two-species Maxwell--Stefan and Busenberg--Travis-type systems and show that for the two-species Maxwell--Stefan equations and the $n$-species thin-film solar cell model, the Boltzmann entropy is unique in some sense (up to a factor and an additive affine function).

The assumption that $\na^2hA$ is symmetric is the same condition that defined so-called extensions of hyperbolic systems $\pa_t u+\pa_x f(u)=0$, with $A$ replaced by the Jacobian $\D f$. It is well known that a strictly convex extension (or entropy) symmetrizes the system \cite{FrLa71}, and that strictly hyperbolic systems of two equations possess large families of extensions \cite{Ser00}. Jenssen and Kogan \cite{JeKo12} showed that the symmetry condition holds if and only if the eigenbasis $(r_i)$ is orthogonal to the inner product defined by $\na^2 h$, and the unknowns become, as in our case, $p_i=r_i^T(\na^2 h)r_i$. These functions satisfy a first-order algebraic-differential system. For $n=2$, it has no algebraic part and coincides with the transport system for $p_\pm$ in Theorem \ref{thm.first}. For $n=3$, we compare our results with the classification of \cite{JeKo12} in Section \ref{sec.jenssen}. 

Our results differ from \cite{JeKo12} in three respects. First, Jenssen and Kogan determine the size of the set of local extensions irrespective of their sign, and they leave the existence of entropies open. We address exactly this question, namely the positivity of $p_\pm$ on (possibly unbounded) global characteristic rectangles. Second, the eigenvalues of a diffusion matrix are not tied to its eigenbasis, as they are for a Jacobian $\D f$. This makes eigenbases relevant that do not matter for conservation laws. Third, we characterize separable convex symmetrizers for any number of components and discuss cross-diffusion systems without symmetric mobility matrix.

The paper is organized as follows. We prove some preliminary results in Section \ref{sec.prep}. The main strategy to construct convex symmetrizers for $n=2$ is detailed in Section \ref{sec.main}. In the case $n=2$, we present sufficient conditions for the existence of convex symmetrizers, provide a complete characterization of separable convex symmetrizers (Section \ref{sec.class}), and construct examples of systems that do not admit an entropy (Section \ref{sec.without}). In Section \ref{sec.exam}, we determine explicit convex symmetrizers and entropy densities for examples from physics and biology, and we present a model for which we go beyond the symmetry condition \eqref{1.hyper}. Finally, in Section \ref{sec.n}, we establish several results for the case $n>2$, including a sufficient condition for the existence of convex symmetrizers and a non-existence criterion.

%%%%%%%%%%%%%%%%%%%%%%%%%%%%%%%%%%%%%%%%%%%%%%%%%%%%%%%%%%%%%%%%%%%

\section{Preliminaries}\label{sec.prep}

We start with some elementary results. Let $\dom\subset\R^2$ be a simply connected open set. Given a positively stable matrix
\begin{align*}
  A(u) = \begin{pmatrix}
  A_{11}(u) & A_{12}(u) \\ A_{21}(u) & A_{22}(u) \end{pmatrix}
  \quad\mbox{for }u=(u_1,u_2)^T\in\dom,
\end{align*}
where $A_{ij}\in C^1(\dom)$, we are looking for a smooth function $h$ with Hessian $\na^2h$ such that $h$ is strictly convex and $\sym(\na^2hA)$ is positive definite in $\dom$, where $\sym M = \frac12(M+M^T)$ denotes the symmetric part of a square matrix $M$. We call a convex function $h\in C^2(\dom)$ an {\em entropy density} if $\sym(\na^2 hA)$ is positive definite and a {\em convex symmetrizer} if $\na^2h(u)A(u)$ is symmetric and positive definite for all $u\in\dom$. Clearly, any convex symmetrizer is an entropy density. Throughout the paper, a function $h\in C^2(\dom)$ is called {\em strictly convex} if its Hessian $\na^2h$ is positive definite in $\dom$. (This is stronger than strict convexity in the usual sense; for instance, $h(u)=u_1^4+u_2^2$ is strictly convex in the usual sense, but its Hessian is singular on $\{u_1=0\}$.)

The following lemma shows that convexity of $\na^2h$ and symmetry of $\na^2h A$ are sufficient to conclude the positive definiteness of $\na^2hA$; see, e.g., \cite[Theorem 6]{ChJu21}.

\begin{lemma}\label{lem.free}
If $A$ is positively stable and $h\in C^2(\dom)$ is strictly convex such that $\na^2hA$ is symmetric in $\dom$, then $\na^2hA$ is positive definite in $\dom$.  
\end{lemma} 

\begin{proof}
We present the short proof for completeness and later reference. Since $H:=\na^2h$ is symmetric, its square root $H^{1/2}$ exists. Set $M=H^{1/2}AH^{-1/2}$. The symmetry of $H$ and $HA$ implies that $M^T  = H^{-1/2}A^TH^{1/2} = H^{-1/2}(A^TH)H^{-1/2} = H^{-1/2}(HA)H^{-1/2} = M$, so $M$ is symmetric too. The matrices $M$ and $A$ are similar, showing that they have the same eigenvalues. As $A$ is positively stable and $M$ is symmetric, $M$ has only positive eigenvalues and is, in fact, positive definite. Moreover, $HA = H^{1/2}(H^{1/2}AH^{-1/2})H^{1/2} = H^{1/2}MH^{1/2}$, implying that $HA$ is congruent to $M$. Since congruence preserves definiteness, the positive definiteness of $M$ implies that $HA$ is positive definite.
\end{proof}

The price for the symmetry requirement is a pointwise obstruction.

\begin{lemma}\label{lem.Delta}
If the function $h$ is strictly convex and $\na^2hA$ is symmetric then
\begin{align*}
  \Delta(u) := (A_{11}-A_{22})^2 + 4A_{12}A_{21} 
  = (\operatorname{tr}A(u))^2-4\det A(u) \ge 0\quad\mbox{for }u\in\dom,
\end{align*}
and $A(u)$ is a scalar multiple of the unit matrix whenever $\Delta(u)=0$. 
\end{lemma}

\begin{proof}
The eigenvalues of $A(u)$ are $\lambda_{\pm}(u) = \frac12(\operatorname{tr}A(u)\pm \sqrt{\Delta(u)})$. The proof of Lemma \ref{lem.free} shows that the eigenvalues of $A(u)$ are real, i.e.\ $\Delta(u)\ge 0$. Let $\Delta(u)=0$. Then $A(u)$ has one double real eigenvalue $\lambda$ (being positive if $A$ is positively stable). We know that $M=H^{1/2}AH^{-1/2}$ is symmetric and has the same eigenvalues as $A$. Then there exists an orthogonal matrix $R$ such that $M=R\operatorname{diag}(\lambda,\lambda)R^T = \lambda I$. This implies that $A = H^{-1/2}MH^{1/2} = \lambda I$.
\end{proof}

Throughout the paper, we impose two conditions (if not stated otherwise):
\begin{align}\label{maincond}
  \Delta(u)>0, \quad A_{12}(u)>0\quad\mbox{for }u\in\dom.
\end{align}
As shown in Lemma \ref{lem.Delta}, the first condition is essential. If $\Delta<0$, convex symmetrizers cannot exist, but entropy densities may still be constructed; see Section \ref{sec.sg} for an example. The second condition in \eqref{maincond} can be replaced by $A_{12}<0$ in $\dom$ (swap the labels $+$ and $-$), or $A_{21}>0$ in $\dom$ (swap the coordinates $u_1$ and $u_2$), or $A_{21}<0$ in $\dom$ (swap both $+$ and $-$ as well as $u_1$ and $u_2$). If $A_{12}$ or $A_{21}$ change sign, convex symmetrizers may not exist globally; see the example given by \eqref{wind} below. 

We show in Remark \ref{rem.locex} below that if $A$ is smooth and positively stable at $u^*$ and $\Delta(u^*)>0$, then a convex symmetrizer exists in a neighborhood of $u^*$. The difficulty is to construct a {\em global} convex symmetrizer. The example \eqref{wind} shows that such a construction may be impossible. 

The positivity of $\Delta$ shows that equation \eqref{1.hyper} is hyperbolic since $\Delta$ is the discriminant associated to \eqref{1.hyper}. Because of $A_{12}>0$, we can normalize the eigenvectors according to
\begin{align*}
  r_\pm = {1\choose m_\pm}, \quad
  m_\pm = \frac{(A_{22}-A_{11})\pm \sqrt{\Delta}}{2A_{12}}.
\end{align*}
The parametrization $r_\pm=(1,m_\pm)^T$ excludes vertical eigendirections, which is why we assume that $A_{12}\neq 0$. The eigenvectors are linearly independent since $\Delta>0$. Thus, they form an eigenbasis in $\R^2$ at every $u\in\dom$. We set $\delta=m_+-m_- = \sqrt{\Delta}/A_{12} > 0$.

%%%%%%%%%%%%%%%%%%%%%%%%%%%%%%%%%%%%%%%%%%%%%%%%%%%%%%%%%%%%%%

\section{Construction of convex symmetrizers ($n=2$)}\label{sec.main}

We reduce the problem of finding a convex symmetrizer to the solution of a linear transport problem along the eigendirections of $A$. To this end, we assume that $A$ is positively stable in $\dom$ and that conditions \eqref{maincond} hold. We introduce the directional derivatives $\D_\pm f$ and the directional curvatures $\kappa_\pm$, respectively, by
\begin{align*}
  \D_\pm f = \pa_1 f + m_\pm\pa_2 f, \quad
  \kappa_\pm = \D_\pm m_\pm,
\end{align*}
recalling that $\pa_i=\pa/\pa u_i$. The curvatures can be interpreted as follows. Let $u_2=y(u_1)$ parametrize an integral curve of $r_\pm$, i.e.\ $y'=m_\pm(u_1,y)$. Then $y'' = \pa_1m_\pm + \pa_2m_\pm\cdot y' = \pa_1m_\pm + m_\pm\pa_2m_\pm = \kappa_\pm$. Thus, $\kappa_+>0$ corresponds to convex $(+)$-characteristic curves, for $\kappa_-<0$ the $(-)$-characteristics are concave curves, and $\kappa_\pm\equiv 0$ means that the family consists of straight lines. 

\subsection{Formulating the Hessian in the eigenbasis}

Let $u\in\dom$ be fixed. We wish to evaluate $\na^2 h=\na^2 h(u)$ in the eigenbasis $(r_\pm)$, which diagonalizes $\na^2h$ if $\na^2hA$ is symmetric (Lemma~\ref{lem.p} below). Let $\pi_\pm = \na h\cdot r_\pm = \pa_1 h + m_\pm\pa_2 h = \D_\pm h$ be the derivative of $h$ along the eigendirection $r_\pm$. 

\begin{lemma}\label{lem.master}
Let $h\in C^2(\dom)$ and $i,j\in\{+,-\}$. Then
\begin{align*}
  r_i^T(\na^2 h)r_j = \D_i\pi_j - (\D_im_j)\pa_2 h\quad\mbox{and}\quad
  \pa_2h = \frac{\pi_+ - \pi_-}{\delta},
\end{align*}
where $\delta=m_+ - m_- = \sqrt{\Delta}/A_{12}>0$. 
\end{lemma} 

\begin{proof}
It follows from the product rule that 
\begin{align*}
  \D_i\pi_j &- (\D_im_j)\pa_2 h 
  = \D_i(\pa_1 h + m_j\pa_2 h) - (\D_im_j)\pa_2 h
  =  \D_i\pa_1 h + m_j\D_i\pa_2 h \\
  &= (\na\pa_1 h)r_i + m_j(\na\pa_2 h)r_i 
  = r_i^T{\na\pa_1 h \choose \na\pa_2 h}\begin{pmatrix} 1 \\ m_j \end{pmatrix}
  = r_i^T(\na^2 h) r_j,
\end{align*}
since the matrix with rows $\na\pa_1 h$, $\na\pa_2 h$ is exactly $\na^2 h$ and $r_j=(1,m_j)^T$. This shows the first identity. The second one follows from $\pi_+ - \pi_- = (\pa_1 h+m_+\pa_2 h)-(\pa_1h+m_-\pa_2 h) = \delta\pa_2 h$.
\end{proof}

The symmetry of $\na^2hA$ is equivalent to the fact that $\na^2h$ diagonalizes in the basis $r_\pm$, and, in this case, the positive definiteness of $\na^2hA$ is equivalent to the positivity of $p_\pm$. 

\begin{lemma}\label{lem.p}
Let $h\in C^2(\dom)$ and set $p_\pm:=r_\pm^T(\na^2h)r_\pm$. Then $p_\pm = \D_\pm\pi_\pm - \kappa_\pm \pa_2 h$. The matrix $\na^2 hA$ is symmetric if and only if 
\begin{align*}
  r_+^T(\na^2 h)r_-=0 \quad \Longleftrightarrow\quad
  \D_+\pi_- = (\D_+m_-)\pa_2 h \quad \Longleftrightarrow\quad
  \D_-\pi_+ = (\D_-m_+)\pa_2 h.
\end{align*}
Let $\na^2hA$ be symmetric. Then:
\begin{enumerate}[label=\rm(\roman*)]
\item $\na^2h$ is positive definite if and only if $p_+>0$ and $p_->0$.
\item If $A$ is positively stable, $\na^2hA$ is positive definite if and only if $p_+>0$ and $p_->0$.
\end{enumerate} 
Moreover, $\det\na^2 h=p_+p_-/\delta^2$. 
\end{lemma}

\begin{proof}
The first statement is an immediate consequence of Lemma \ref{lem.master}, since $p_\pm = \D_\pm\pi_\pm - (\D_\pm m_\pm)\pa_2h = \D_\pm\pi_\pm - \kappa_\pm\pa_2 h$. For the second statement, if $\na^2hA$ is symmetric, we conclude from $\lambda_+\neq\lambda_-$ and
\begin{align*}
  0 = r_+^T(A^T\na^2 h - \na^2hA)r_- 
  = (\lambda_+-\lambda_-)r_+^T(\na^2 h)r_-
\end{align*}
that $r_+^T(\na^2 h)r_-=0$. On the other hand, if $r_+^T(\na^2 h)r_-=0$ then $r_+^T Br_-=0$, where $B=A^T\na^2 h - \na^2hA$. Moreover, $r_-^TBr_+ = - r_+^TBr_-=0$ and $r_+^TBr_+ = r_-^TBr_- = 0$. The four numbers $r_i^TBr_j$ determine the matrix $B$ completely, giving $B=0$. This proves the symmetry of $\na^2hA$. The reformulations follow from Lemma \ref{lem.master} with $(i,j)=(+,-)$ and $(i,j)=(-,+)$. For the third statement, since $r_+^T(\na^2h)r_-=0$, the matrix $\na^2h$ becomes $\operatorname{diag}(p_+,p_-)$ in the basis $(r_+,r_-)$. In particular, $\na^2h$ is positive definite if and only if $p_\pm>0$, which proves (i). 

For (ii), let $A$ be positively stable. By (i), it suffices to show that $\na^2hA$ is positive definite if and only if $\na^2h$ is positive definite. If $\na^2h$ is positive definite, then $\na^2hA$ is positive definite by Lemma~\ref{lem.free}. Conversely, let $H:=\na^2h$ be such that $HA$ is symmetric and positive definite. By \cite[Problem 9, Section 2.2]{HoJo91},
\begin{align*}
  H = \int_0^\infty e^{-tA^T}(A^TH+HA)e^{-tA}\dd t
  = 2\int_0^\infty e^{-tA^T}HAe^{-tA}\dd t,
\end{align*}
which is positive definite, since $HA$ does so. Finally, writing $R=(r_+,r_-)\in\R^{2\times 2}$, we have $\det R=m_- - m_+ = -\delta$, $\operatorname{diag}(p_+,p_-) = R^T(\na^2 h)R$, and hence
\begin{align*}
  p_+p_- = \det\operatorname{diag}(p_+,p_-) = \det(R^T(\na^2 h)R)
  = (\det R)^2\det(\na^2 h) = \delta^2\det(\na^2h),
\end{align*}
which finishes the proof.
\end{proof}

Lemma \ref{lem.p} says that $\na^2hA$ is symmetric when the derivative of $h$ along one characteristic direction is constrained as one moves
in the other characteristic direction. This is a statement about first derivatives only and therefore, the problem is in fact a first-order system. The following lemma shows that the Hessian $\na^2h$ can be reconstructed from the values $p_\pm$. We show this in a more general context.

\begin{lemma}[Reconstruction]\label{lem.recon}
Let $p_+$, $p_-\in\R$ and $H=(H_{ij})\in\R^{2\times 2}$ be the unique symmetric matrix with $r_+^THr_+ = p_+$, $r_-^THr_- = p_-$, $r_+^THr_- = 0$. Then
\begin{align*}
  H_{11} = \frac{p_+m_-^2+p_-m_+^2}{\delta^2}, \quad
  H_{12} = -\frac{p_+m_-+p_-m_+}{\delta^2}, \quad
  H_{22} = \frac{p_+ + p_-}{\delta^2}.
\end{align*}
\end{lemma}

\begin{proof}
The three conditions on $H$ imply that $R^THR=\operatorname{diag}(p_+,p_-)$, where $R=(r_+,r_-)\in\R^{2\times 2}$. Thus, $H = R^{-T}\operatorname{diag}(p_+,p_-)R^{-1}$. The inverse $R^{-1}$ is given by the dual basis $(\rho_+,\rho_-)$, where the row vectors $\rho_\pm$ are defined by $\rho_\pm r_\pm=1$, $\rho_+ r_- = \rho_- r_+ = 0$. More explicitly, the dual basis reads as $\rho_+=\delta^{-1}(-m_-,1)$, $\rho_-=\delta^{-1}(m_+,-1)$. Then
\begin{align}\label{H}
  H = p_+(\rho_+^T \rho_+) + p_-(\rho_-^T \rho_-)
  = \frac{p_+}{\delta^2}\begin{pmatrix}
  m_-^2 & -m_- \\ -m_- & 1 \end{pmatrix}
  + \frac{p_-}{\delta^2}\begin{pmatrix}
  m_+^2 & -m_+ \\ -m_+ & 1 \end{pmatrix},
\end{align}
which completes the proof.
\end{proof}

%%%%%%%%%%%%%%%%%%

\subsection{First-order system}

The aim of this section is to derive a first-order system for $(p_+,p_-)$, which allows us to reconstruct the Hessian $\na^2 h$ via Lemma \ref{lem.recon}. We first show that $h\mapsto(\pi_+,\pi_-)$ is a bijection between solutions to the hyperbolic equation \eqref{1.hyper} and the first-order system \eqref{pi} below for $\pi_\pm$. 

\begin{lemma}[First-order formulation]
The map $h\mapsto (\pi_+,\pi_-)$ is a bijection, up to additive constants, between solutions $h\in C^2(\dom)$ to \eqref{1.hyper} and solutions $(\pi_+,\pi_-)\in C^1(\dom)^2$ to the first-order system
\begin{align}\label{pi}
  \D_+\pi_- = \frac{\pi_+-\pi_-}{\delta}\D_+m_-, \quad
  \D_-\pi_+ = \frac{\pi_+-\pi_-}{\delta}\D_-m_+.
\end{align}
\end{lemma}

\begin{proof}
Let $h\in C^2(\dom)$ be a solution to \eqref{1.hyper}. Then Lemmas \ref{lem.p} and \ref{lem.master} show that 
\begin{align*}
  \D_+\pi_- = (\D_+m_-)\pa_2 h = (\D_+m_-)\frac{\pi_+-\pi_-}{\delta}, \quad
  \D_-\pi_+ = (\D_-m_+)\frac{\pi_+-\pi_-}{\delta},
\end{align*}
which is \eqref{pi}. Conversely, if $(\pi_+,\pi_-)\in C^1(\dom)^2$ is given, we define $g_2:=(\pi_+-\pi_-)/\delta$ and $g_1:=\pi_+ - m_+g_2$. Then $\pi_\pm = g_1+m_\pm g_2$. Computing exactly as in the proof of Lemma \ref{lem.master} (which used only the product rule, not the existence of $h$), we obtain
\begin{align*}
  \D_i\pi_j - (\D_im_j)g_2 = r_i^TG^Tr_j, \quad\mbox{where }
  G = {\na g_1 \choose \na g_2}\mbox{ or }G_{ij} = \pa_jg_i
\end{align*}
for $i$, $j\in\{+,-\}$. Equations \eqref{pi} imply that $r_+^TG^Tr_- = r_-^TG^Tr_+ = 0$ and hence $r_+^T(G-G^T)r_-=0$. Thus, since $G-G^T$ is antisymmetric, this matrix is a multiple of $J=\left(\begin{smallmatrix} 0 & -1 \\ 1 & 0 \end{smallmatrix}\right)$. We infer from $r_+^TJr_- = m_+-m_-=\delta\neq 0$ that the multiple must be zero, so $G-G^T=0$. Hence, $\pa_2g_1=\pa_1g_2$ in $\dom$. Since $\dom$ is simply connected by assumption, Poincar\'e's lemma implies the existence of a function $h$ such that $\na h=(g_1,g_2)$, being unique up to a constant. Thus, 
\begin{align*}
  0 &= r_+^TGr_- = r_+^T(\na^2 h)r_- 
  = \begin{pmatrix} 1 \\ m_+ \end{pmatrix}^T
  \begin{pmatrix} \pa_{11}^2 h & \pa_{12}^2 h \\ 
  \pa_{12}^2 h & \pa_{22}^2 h \end{pmatrix}
  \begin{pmatrix} 1 \\ m_- \end{pmatrix} \\
  &= \pa_{11}^2 h + (m_+ + m_-)\pa_{12}^2 h + m_+m_-\pa_{22}^2 h
  = \pa_{11}^2 h + \frac{A_{22}-A_{11}}{A_{12}} \pa_{12}^2 h
  - \frac{A_{21}}{A_{12}}\pa_{22}^2 h,
\end{align*}
which equals \eqref{1.hyper} after multiplication by $A_{12}$.
\end{proof}

We prove that the map $h\mapsto(p_+,p_-)$ is a bijection. For this, we first compute the derivative of the basis vectors. Recalling that $r_\pm = (1,m_\pm)^T$ has a constant first component, we find for $i$, $j\in\{+,-\}$ that
\begin{align}\label{Dr}
  r_+ - r_- = \delta\begin{pmatrix} 0 \\ 1 \end{pmatrix}, \quad
  \D_ir_j = \begin{pmatrix} 0 \\ \D_im_j \end{pmatrix}
  = \frac{\D_im_j}{\delta}(r_+ - r_-). 
\end{align}

\begin{theorem}[First-order formulation for $(p_+,p_-)$]\label{thm.first}
The map $h\mapsto (p_+,p_-)$ is a bijection, modulo addition of affine functions, between solutions $h\in C^2(\dom)$ to \eqref{1.hyper} and solutions $(p_+,p_-)\in C^0(\dom)^2$ with $\D_\pm p_\mp\in C^0(\dom)$ to the first-order system
\begin{equation}\label{star}
\begin{aligned}
  \D_- p_+ &= \frac{2\D_-m_+ - \D_+m_-}{\delta}p_+ 
  + \frac{\kappa_+}{\delta}p_-, \\
  \D_+ p_- &= -\frac{2\D_+m_- - \D_-m_+}{\delta}p_- 
  - \frac{\kappa_-}{\delta}p_+,
\end{aligned}
\end{equation}
recalling that $\kappa_\pm=\D_\pm m_\pm$. Moreover, $h$ is a {\em convex} symmetrizer if and only if $p_+>0$ and $p_->0$. 
\end{theorem}

\begin{proof}
To simplify, we suppose that $h\in C^3(\dom)$ and $p_\pm\in C^1(\dom)$. The arguments can be extended to functions $h\in C^2(\dom)$ and $p_\pm\in C^0(\dom)$ with continuous characteristic derivatives (see Step 3 of the proof). 

{\em Step 1:} We differentiate $r_j^T(\na^2 h)r_k$ along $r_i$, where $i$, $j$, $k\in\{+,-\}$, using the product rule:
\begin{align}\label{2.aux}
  \D_i(r_j^T(\na^2 h)r_k) = K_{ijk} + (\D_ir_j)^T(\na^2 h)r_k
  + r_j^T(\na^2 h)\D_ir_k, 
\end{align}
where $K_{ijk} = r_j^T(\D_i\na^2 h)r_k = \sum_{a,b,c=1}^2r_{i,a} r_{j,b} r_{k,c}\pa_{abc}^3 h$ is symmetric under permutations of $i$, $j$, $k$ (here, we need $h\in C^3(\dom)$ and write $r_{i}=(r_{i,1},r_{i,2})^T$). By \eqref{Dr}, the second term on the right-hand side of \eqref{2.aux} can be written as
\begin{align*}
  (\D_ir_j)^T(\na^2 h)r_k = \frac{\D_im_j}{\delta}
  \big(r_+^T(\na^2 h)r_k - r_-^T(\na^2 h)r_k\big),
\end{align*}
which, using the symmetry of $\na^2h$, leads to
\begin{align}\label{Dr2}
  \D_i(r_j^T(\na^2 h)r_k) &= K_{ijk} 
  + \frac{\D_im_j}{\delta}
  \big(r_+^T(\na^2 h)r_k - r_-^T(\na^2 h)r_k\big) \\
  &\phantom{xx}+ \frac{\D_im_k}{\delta}
  \big(r_+^T(\na^2 h)r_j - r_-^T(\na^2 h)r_j\big). \nonumber 
\end{align}
Choosing $i=-$, $j=k=+$ and taking into account that $r_+^T(\na^2 h)r_+ = p_+$ and $r_-^T(\na^2 h)r_+=0$, this becomes
\begin{align}\label{2.Dp}
  \D_-p_+ &= \D_-(r_+^T(\na^2 h)r_+) = K_{-++} + 2\frac{\D_-m_+}{\delta}
  \big(r_+^T(\na^2 h)r_+ - r_-^T(\na^2 h)r_+\big) \\
  &= K_{-++} + 2\frac{\D_-m_+}{\delta}p_+. \nonumber 
\end{align}
To determine $K_{-++}=K_{++-}$, we differentiate $r_+^T(\na^2 h)r_-=0$ and use \eqref{Dr2} with $i=j=+$, $k=-$ as well as $r_+^T(\na^2 h)r_-=0$:
\begin{align*}
  0 &= \D_+\big(r_+^T(\na^2 h)r_-\big)
  = K_{++-} + \frac{\kappa_+}{\delta}\big(r_+^T(\na^2 h)r_-
  - r_-^T(\na^2 h)r_-\big) \\
  &\phantom{xx}+ \frac{\D_+m_-}{\delta}
  \big(r_+^T(\na^2 h)r_+ - r_-^T(\na^2 h)r_+\big) 
  = K_{++-} - \frac{\kappa_+}{\delta}p_- + \frac{\D_+m_-}{\delta}p_+.
\end{align*}
It follows that $ K_{-++} = K_{++-} = \delta^{-1}\kappa_+p_- - \delta^{-1}\D_+m_-p_+$. Substituting this expression into \eqref{2.Dp} gives the first line of \eqref{star}. The second line is obtained by exchanging $+$ and $-$, which replaces $\delta$ by $-\delta$ and swaps $p_+$ and $p_-$. 

{\em Step 2:} Let $(p_+,p_-)$ solve \eqref{star}. By Lemma \ref{lem.recon}, there exists a unique symmetric matrix $H$ such that $r_+^THr_+ = p_+$, $r_-^THr_- = p_-$, and $r_+^THr_-=0$. As $H$ is determined from $p_\pm$, the regularity $p_\pm\in C^1(\dom)$ implies that $H\in C^1(\dom)$. We use the following identities, which hold for {\em every} symmetric matrix field $T\in C^1(\dom)$, not necessarily a Hessian. Set $q_\pm:=r_\pm^TTr_\pm$, $t:=r_+^TTr_-$, $c=(c_1,c_2)^T$ with $c_i:=\pa_1 T_{i2}-\pa_2 T_{i1}$ for $i=1,2$, and $\beta_+:=(2\D_-m_+-\D_+m_-)/\delta$,  $\beta_-:=(2\D_+m_--\D_-m_+)/\delta$. Using the product rule and identities \eqref{Dr}, a computation shows that
\begin{equation}\label{ident}
\begin{aligned}
  \D_-q_+ - \D_+t &= \beta_+(q_+-t) + \frac{\kappa_+}{\delta}(q_--t) 
  + \delta(r_+\cdot c), \\
  \D_+q_- - \D_-t &= -\beta_-(q_--t) - \frac{\kappa_-}{\delta}(q_+-t) 
  - \delta(r_-\cdot c). 
\end{aligned}
\end{equation}
%Indeed, denoting by $\D_\pm T$ the entrywise derivative of $T$, the product rule gives
%\begin{align*}
%  \D_-q_+ - \D_+t = r_+^T\big[(\D_-T)r_+ - (\D_+T)r_-\big]
%  + 2(\D_-r_+)^TTr_+ - (\D_+r_+)^TTr_- - r_+^TT\D_+r_-.
%\end{align*}
%Since $\D_\pm T=\pa_1T+m_\pm\pa_2T$, $r_+-r_-=\delta(0,1)^T$, and $m_-r_+-m_+r_-=-\delta(1,0)^T$, we have $(\D_-T)r_+-(\D_+T)r_- = \delta(\pa_1Te_2-\pa_2Te_1) = \delta c$, so the first term equals $\delta(r_+\cdot c)$. By \eqref{Dr}, the remaining three terms equal $\frac{2\D_-m_+}{\delta}(q_+-t) - \frac{\kappa_+}{\delta}(t-q_-) - \frac{\D_+m_-}{\delta}(q_+-t)$, which gives the first identity. The second one follows in the same way.
We apply these identities to $T=H$, for which $q_\pm=p_\pm$ and $t=0$. In view of \eqref{star} and $\delta\neq 0$, they give $r_+\cdot c=0$ and $r_-\cdot c=0$. Since $(r_+,r_-)$ is a basis, this shows that $c=0$ and $\pa_2 H_{1i} = \pa_1 H_{2i}$. The symmetry of the $2\times 2$ matrix $H$ then yields the total symmetry of $\pa_i H_{jk}$. 

It remains to prove that $H$ is a Hessian. To this end, let $k\in\{1,2\}$ and set $F^{(k)} = (H_{1k},H_{2k})$. By total symmetry, we have $\pa_2 F^{(k)}_1 = \pa_2 H_{1k} = \pa_1 H_{2k} = \pa_1 F^{(k)}_2$. Thus, by the Poincar\'e lemma, there exists a function $g_k$ such that $\pa_j g_k = H_{jk}$ for $j=1,2$. Because of $\pa_2g_1 = H_{21} = H_{12} = \pa_1 g_2$, we can apply the Poincar\'e lemma a second time, and there exists a function $h$ such that $g_k=\pa_k h$. This shows that $\pa_{jk}^2h=\pa_j g_k=H_{jk}$, i.e.\ $\na^2 h = H$. Since $\na^2 h$ determines $h$ up to an affine function, the correspondence is bijective modulo affine functions. Finally, the last statement follows from the last statement in Lemma \ref{lem.p}. 

{\em Step 3:} We outline the arguments for relaxing the regularity on $h$ and $p_\pm$. Identities \eqref{ident} are linear and of first order in $T$, with smooth coefficients. Therefore, approximating $T$ by mollified fields, they hold in the sense of distributions for every {\em continuous} symmetric matrix field $T$. 
%(Note that the Hessian is mollified, not the equation \eqref{1.hyper}, which is not preserved under mollification.) 
This allows us to treat $h\in C^2(\dom)$. Indeed, if $h\in C^2(\dom)$ solves \eqref{1.hyper}, then $T=\na^2h$ is continuous, $t=0$, and $c=0$ in the sense of distributions, since distributional derivatives commute. Hence, $\D_-p_+$ equals a continuous function in the sense of distributions. 
%In the characteristic coordinates of Lemma~\ref{lem.char} below, $\D_-=A_-\pa_{\alpha_-}$, and a continuous function whose distributional derivative $\pa_{\alpha_-}$ is continuous is classically differentiable in $\alpha_-$. 
Thus, the first equation in \eqref{star} holds pointwise (and similarly for the second one). Conversely, let $p_\pm\in C^0(\dom)$ have continuous derivatives $\D_-p_+$, $\D_+p_-$ and solve \eqref{star}. Then the matrix $H$ of Lemma~\ref{lem.recon} is continuous, $t=0$, and the identities give $c=0$ in the sense of distributions. The Poincar\'e lemma also holds for continuous vector fields with vanishing distributional curl, and we obtain, as in Step~2, a function $h\in C^2(\dom)$ with $\na^2h=H$.
\end{proof}

Equations \eqref{star} determine $\D_-p_+$ and $\D_+p_-$ only. So, $p_+$ is transported across the $(+)$-characteristic and $p_-$ across the $(-)$-characteristic, while $\D_+p_+$ and $\D_-p_-$ are unconstrained. Hence, prescribing $p_+>0$ on one $(+)$-characteristic $C_+$ and $p_->0$ on one $(-)$-characteristic $C_-$, the function $h$ is determined up to an affine function, on the rectangle spanned by $C_+$ and $C_-$. 

%%%%%%%%%%%%%%%%%%%%%%%

\subsection{Normal form}

We can remove the diagonal coefficients of problem \eqref{star} by a change of variables and an integrating factor. We say that a function $\alpha$ is a first integral of $r_\pm$ if $\D_\pm\alpha=0$, i.e.\ if $\alpha$ is constant along every $(\pm)$-characteristic curve. 

\begin{lemma}[Characteristic coordinates]\label{lem.char}
Let assumption \eqref{maincond} hold. Locally (and globally on $\dom$ whenever the transport problems below are globally solvable), there are functions $\alpha_+$, $\alpha_-$ such that
\begin{align*}
  \D_+\alpha_- = \D_-\alpha_+ = 0, \quad \D_+\alpha_+ =:A_+>0, \quad
  \D_-\alpha_- =:A_->0.
\end{align*}
Moreover, recalling that $(\rho_+,\rho_-)$ is the dual basis of $(r_+,r_-)$ (see the proof of Lemma \ref{lem.recon}),
\begin{enumerate}[label=\rm(\roman*)]
\item $\na\alpha_\pm = A_\pm\rho_\pm$; 
%\item[\rm (2)] the $(+)$-characteristics are the level sets of $\alpha_-$ and the $(-)$-characteristics are the level sets of $\alpha_+$;
\item $\D_+ = A_+\pa_{\alpha_+}$ and $\D_- = A_-\pa_{\alpha_-}$. 
\end{enumerate}
\end{lemma}

Property (i) means that the Jacobian of $u\mapsto(\alpha_+,\alpha_-)$ is invertible. Thus, $(\alpha_+,\alpha_-)$ is a local coordinate system. 

\begin{proof}
The existence of local first integrals with $A_\pm\neq 0$ is a consequence of the flow-box theorem, which states that in the neighborhood of a regular point (where the vector field does not vanish), a smooth vector field can be locally transformed into a constant vector field. Replacing $\alpha_\pm$ by $-\alpha_\pm$ (if necessary) makes $A_\pm>0$.

(1) Let $\eta$ be a row vector. It is determined by $\eta = (\eta r_+)\rho_+ + (\eta r_-)\rho_-$. We apply this to $\eta=\na\alpha_+$. Since $(\na\alpha_+) r_+ = \D_+\alpha_+ = A_+$ and $(\na\alpha_+) r_- = \D_-\alpha_+=0$, we obtain $\na\alpha_+ = A_+\rho_+$. Similarly, we have $\na\alpha_-=A_-\rho_-$. This means that the Jacobian of $(\alpha_+,\alpha_-)$ has the rows $A_+\rho_+$ and $A_-\rho_-$, which are linearly independent because $\rho_\pm$ is a basis and $A_\pm\neq 0$. By the inverse function theorem, $(\alpha_+,\alpha_-)$ are coordinates. 

%(2) The level set $\{\alpha_-=\mbox{const.}\}$ is a curve whose tangent is annihilated by $\na\alpha_-=A_-\rho_-$. As $\rho_-r_+=0$, the tangent is in fact $r_+$. So, the level sets of $\alpha_-$ are the $(+)$-characteristics. A similar statement holds for the level sets of $\alpha_+$.

(2) For any function $f\in C^1$, the chain rule gives $\D_\pm f = (\D_\pm\alpha_+)\pa_{\alpha_+}f + (\D_\pm \alpha_-)\pa_{\alpha_-}f$. Because of $\D_-\alpha_+=\D_+\alpha_-=0$, this becomes $\D_+ f= A_+\pa_{\alpha_+}f$ and $\D_-f = A_-\pa_{\alpha_-}f$. 
\end{proof} 

In the $(\alpha_+,\alpha_-)$ coordinates, the rectangle $R=[a_0,a_1]\times[b_0,b_1]$ is bounded by two $(+)$-characteristics $(\alpha_-=b_0,b_1)$ and two $(-)$-characteristics $(\alpha_+=a_0,a_1)$. Its preimage in $\dom$ is the curvilinear quadrilateral spanned by these four curves, called the characteristic rectangle. Prescribing data on one side of each pair is the Goursat problem \cite[Chap.~6]{RuRu94}.

The characteristic rectangle $R$ may be unbounded: We call a set $R\subset\dom$ a {\em characteristic rectangle} if $u\mapsto(\alpha_+(u),\alpha-(u))$ maps $R$ bijectively onto $I\times J$ of intervals $I$, $J\subset\R$, and we identify $R=I\times J$. If $I$ or $J$ are unbounded, the corresponding side ``lies at infinity'', and Goursat data on such a side are understood as limits.   

\begin{proposition}[Normal form]\label{prop.nf}
On a characteristic rectangle $R=[a_0,a_1]\times[b_0,b_1]$, there exist positive functions $a_\pm$ and functions $K_\pm$ such that system \eqref{star} becomes
\begin{align}\label{star2}
  \pa_{\alpha_-}P_+ = K_+P_-, \quad \pa_{\alpha_+}P_- = K_-P_+, 
\end{align}
where $P_\pm = p_\pm/a_\pm$ and $\operatorname{sign}K_\pm = \pm\operatorname{sign}\kappa_\pm$. Moreover, $h$ is a convex symmetrizer if and only if $P_+>0$ and $P_->0$, and 
%the product $K_+K_-$ is independent of the choice of $a_\pm$, with 
$K_+K_- = -\kappa_+\kappa_-/(\delta^2 A_+A_-)$.
\end{proposition}

\begin{proof}
First, we rewrite \eqref{star} in the coordinates $(\alpha_+,\alpha_-)$. By Lemma \ref{lem.char} (2), we have $\D_-p_+ = A_-\pa_{\alpha_-}p_+$ and $\D_+p_- = A_+\pa_{\alpha_+}p_-$. We divide the first line of \eqref{star} by $A_->0$ and the second one by $A_+>0$:
\begin{align}\label{pp}
  & \pa_{\alpha_-}p_+ = L_{11}p_+ + L_{12}p_-, \quad
  \pa_{\alpha_+}p_- = L_{21}p_- + L_{22}p_+, \quad\mbox{where} \\
  & L_{11} = \frac{2\D_-m_+-\D_+m_-}{\delta A_-}, \
  L_{12} = \frac{\kappa_+}{\delta A_-}, \
  L_{21} = -\frac{2\D_+m_- - \D_-m_+}{\delta A_+}, \
  L_{22} = -\frac{\kappa_-}{\delta A_+}. \nonumber 
\end{align}
Next, we substitute $p_\pm = a_\pm P_\pm$ in \eqref{pp}. By the product rule, $\pa_{\alpha_-}(a_+P_+) = a_+\pa_{\alpha_-}P_+ + (\pa_{\alpha_-}a_+)P_+$. Thus, the first equation in \eqref{pp} becomes, after division by $a_+>0$,
\begin{align*}
  \pa_{\alpha_-}P_+ &= \bigg(L_{11} 
  - \frac{\pa_{\alpha_-}a_+}{a_+}\bigg)P_+ + \frac{L_{12}a_-}{a_+}P_-,
  \quad\mbox{and similarly,} \\
  \pa_{\alpha_+}P_- &= \bigg(L_{21} 
  - \frac{\pa_{\alpha_+}a_-}{a_-}\bigg)P_- + \frac{L_{22}a_+}{a_-}P_+.
\end{align*}

We aim to eliminate the brackets, i.e., we are seeking $a_\pm$ such that $L_{11}=\pa_{\alpha_-}a_+/a_+=\pa_{\alpha_-}\log a_+$ and $L_{21}$ $= \pa_{\alpha_+}a_-/a_- = \pa_{\alpha_+}\log a_-$. These are ordinary differential equations in one variable, whose solutions are given by
\begin{align*}%\label{apm}
  a_+(\alpha_+,\alpha_-) = \exp\int_{b_0}^{\alpha_-}
  L_{11}(\alpha_+,s)\dd s, \quad
  a_-(\alpha_+,\alpha_-) = \exp\int_{a_0}^{\alpha_+}
  L_{21}(s,\alpha_-)\dd s.
\end{align*}
With this choice, system \eqref{star} becomes \eqref{star2}, where
\begin{align*}
  K_+ = \frac{L_{12}a_-}{a_+}, \quad K_- = \frac{L_{22}a_+}{a_-}.
\end{align*}

Since $a_\pm>0$ and $\delta>0$, $A_\pm>0$, the signs of $K_+$ and $K_-$ are those of $L_{12}$ and $L_{22}$, i.e.\ those of $\kappa_+$ and $-\kappa_-$, respectively. The factors $a_\pm$ cancel in the product, $K_+K_- = L_{12}L_{22} = -\kappa_+\kappa_-/(\delta^2 A_+A_-)$. 
%which is therefore independent of the choice of $a_0$ and $b_0$. 
Finally, $P_\pm>0$ if and only if $p_\pm>0$ (since $a_\pm>0$ and $p_\pm = a_\pm P_\pm$), and the last claim is the last statement in Lemma \ref{lem.p}. 
\end{proof}

We see from the proof that $P_+$ changes only as one moves in $\alpha_-$ (across $(+)$-characteristics), and it changes by $K_+$ times $P_-$; similarly for $P_-$. The only question is whether the two coupling terms can drive $P_+$ or $P_-$ through zero, and that depends solely on the sizes and signs of $K_+$ and $K_-$. 

\begin{lemma}[Well-posedness]\label{lem.well}
Let $K_\pm$ be continuous on the compact rectangle $R = [a_0,a_1]\times[b_0,b_1]$. Let $\phi_+\in C^0([a_0,a_1])$ (the values of $P_+$ on the edge $\alpha_-=b_0$) and $\phi_-\in C^0([b_0,b_1])$ (the values of $P_-$ on the edge $\alpha_+=a_0$) be given. Then the integral equations
\begin{equation}\label{inteq}
\begin{aligned}
  P_+(\alpha_+,\alpha_-) &= \phi_+(\alpha_+)
  + \int_{b_0}^{\alpha_-} K_+(\alpha_+,s)P_-(\alpha_+,s)\dd s, \\
  P_-(\alpha_+,\alpha_-) &= \phi_-(\alpha_-)
  + \int_{a_0}^{\alpha_+} K_-(s,\alpha_-)P_+(s,\alpha_-)\dd s
\end{aligned}
\end{equation}
have a unique continuous solution on $R$, which solves the Goursat problem \eqref{star2} on $R$ with initial data $P_+(\cdot,b_0)=\phi_+$, $P_-(a_0,\cdot)=\phi_-$. If $\phi_\pm\in C^k$ and $K_\pm\in C^k(R)$ for some $k\ge 1$, then $P_\pm\in C^k(R)$. 
\end{lemma}

\begin{proof}
We introduce on $X:=C^0(R)^2$ the norm 
\begin{align*}
  \|(P_+,P_-)\|_\lambda := \sup_{(\alpha_+,\alpha_-)\in R} e^{-\lambda(\alpha_+ + \alpha_-)}
  \big(|P_+(\alpha_+,\alpha_-)|+|P_-(\alpha_+,\alpha_-)|\big), 
\end{align*}
where $\lambda>0$. Clearly, $(X,\|\cdot\|_\lambda)$ is a Banach space. Define the map $F:X\to X$ by $F=(F_+,F_-)$, where 
\begin{align*}
  F_+(P_+,P_-)(\alpha_+,\alpha_-) 
  &= \phi_+(\alpha_+) + \int_{b_0}^{\alpha_-} 
  K_+(\alpha_+,s)P_-(\alpha_+,s)\dd s, \\ 
  F_-(P_+,P_-)(\alpha_+,\alpha_-) 
  &= \phi_-(\alpha_-) + \int_{a_0}^{\alpha_+} 
  K_-(s,\alpha_-)P_+(s,\alpha_-)\dd s.
\end{align*}
We claim that for sufficiently large $\lambda>0$, $F$ is a contraction on $X$. Indeed, 
\begin{align*}
  |F_+(P_+,P_-)-F_+(P_+^*,P_-^*)|(\alpha_+,\alpha_-)
  &\le \sup_R|K_+|\|(P_+-P_+^*,P_--P_-^*)\|_\lambda \int_{b_0}^{\alpha_-} 
  e^{\lambda(\alpha_++s)}\dd s \\
  &\le \frac{e^{\lambda(\alpha_++\alpha_-)}}{\lambda}\sup_R|K_+|
  \|(P_+-P_+^*,P_--P_-^*)\|_\lambda,
\end{align*}
and similarly for $F_-$. This gives
\begin{align*}
  \|F(P_+,P_-)-F(P_+^*,P_-^*)\|_\lambda
  \le \frac{C}{\lambda}\|(P_+-P_+^*,P_--P_-^*)\|_\lambda
\end{align*}
for $C=\sup_R|K_+|+\sup_R|K_-|$. Choosing $\lambda>C$, $F$ becomes a contraction on $X$, and we can apply the Banach fixed-point theorem. The fixed point solves \eqref{star2} with the given data. For the regularity, let $\phi_\pm\in C^1$ and $K_\pm\in C^1(R)$. Differentiating \eqref{inteq} and using $\pa_{\alpha_+}P_-=K_-P_+$, $\pa_{\alpha_-}P_+=K_+P_-$, we find that
\begin{align*}
  \pa_{\alpha_+}P_+ &= \phi_+' + \int_{b_0}^{\alpha_-}\big(\pa_{\alpha_+}K_+P_- + K_+K_-P_+\big)(\alpha_+,s)\dd s, \\
  \pa_{\alpha_-}P_- &= \phi_-' + \int_{a_0}^{\alpha_+}\big(\pa_{\alpha_-}K_-P_+ + K_-K_+P_-\big)(s,\alpha_-)\dd s.
\end{align*}
The right-hand sides are continuous, so $P_\pm\in C^1(R)$. The case $k\ge 2$ follows by induction.
\end{proof}

Lemma \ref{lem.well} prescribes data on the edges through the corner $(a_0,b_0)$. Data on the edges through another corner are reduced to this case by replacing $\alpha_+$ by $-\alpha_+$, or $\alpha_-$ by $-\alpha_-$. Such a reflection changes the sign of $K_-$, or of $K_+$, but not their bounds, so the lemma applies unchanged.

\begin{remark}[Local existence]\label{rem.locex}
The local existence of a convex symmetrizer is automatic: Let $A$ be smooth and positively stable near $u^*\in\dom$ and let \eqref{maincond} hold at that point. Consider a neighborhood of $u^*$ in which a smooth eigenbasis and characteristic coordinates are defined, and solve \eqref{star2} with data $P_\pm\equiv 1$ on the two characteristics through $u^*$ (see Lemma \ref{lem.well}). Then $P_\pm>0$ on a possibly smaller rectangle. We conclude from Proposition \ref{prop.nf} the existence of a convex symmetrizer in a neighborhood of $u^*$, constructed from Lemma \ref{lem.recon}. 
\end{remark}

More generally, let $\dom$ be a simply connected domain. We have proved the existence of a solution $h$ to the symmetrizer equation \eqref{1.hyper} on any compact characteristic rectangle. If the solution to the transport system stays positive, $h$ is strictly convex and consequently a convex symmetrizer. The convex symmetrizer exists on $\dom$ if the pair of data curves reaches all of $\dom$ (this is important if $\dom$ is unbounded). These two properties (positivity and global characteristic rectangles) have to come from additional results. We discuss some aspects in Section \ref{sec.fail}. Sufficient conditions for the global construction are given in Section \ref{sec.suff}.

\subsection{Example}\label{sec.ex2}

We apply the strategy of the previous subsection to the positively stable matrix
\begin{align*}%\label{ex.full}
  A(u) = \begin{pmatrix} 1+u_1 & 1 \\ u_1 & 2 \end{pmatrix}
  \quad\mbox{for }u\in\dom=(0,\infty)^2.
\end{align*}
We have $\Delta(u)=(1+u_1)^2>0$ and $A_{12}=1>0$, so the main assumptions are fulfilled. 

{\em Step 1: First-order system.}
The slopes $m_\pm$ are the solutions to $m^2+(u_1-1)m-u_1 = (m-1)(m+u_1)=0$. Hence,
\begin{align*}
  m_+ = 1, \quad m_- = -u_1, \quad \delta= 1+u_1, \quad
  \D_+=\pa_1+\pa_2, \quad \D_-=\pa_1-u_1\pa_2.
\end{align*}
The $(+)$-characteristics are the integral curves $u_2=y(u_1)$ determined from $y'(u_1)=m_+=1$ giving the lines $u_2=u_1+\mbox{const}$. The $(-)$-characteristics are computed from $y'(u_1)=m_-=-u_1$; they are the downward parabolas $u_2=\mbox{const.}-u_1^2/2$. Furthermore,
\begin{align*}
  \kappa_+ = \D_+m_+ = \pa_1m_+ + m_+\pa_2m_+ = 0, \quad
  \kappa_- = \D_-m_- = \pa_1m_- + m_-\pa_2m_- = -1
\end{align*} 
and $\D_-m_+=0$, $\D_+m_-=-1$. Equations \eqref{star} become
\begin{align}\label{ex.p}
  \D_-p_+ = \frac{p_+}{1+u_1}, \quad 
  \D_+p_- = \frac{2p_- + p_+}{1+u_1}.
\end{align}
This system is triangular; it is sufficient to solve the first equation and then substitute the solution into the second equation. 

{\em Step 2: Solution of the first-order system.}
The characteristic coordinates from Lemma \ref{lem.char} are determined from $u_2+u_1^2/2 = \mbox{const.}$, $u_2-u_1=\mbox{const.}$, so 
\begin{align*}
  \alpha_+ = u_2 + \tfrac12 u_1^2, \quad \alpha_- = u_1-u_2.
\end{align*} 
The numbers $A_+ = \D_+\alpha_+=1+u_1$, $A_- = \D_-\alpha_- = 1+u_1$ are both positive. Since $\D_-(1+u_1)=1$, the function $1+u_1$ solves the first equation of \eqref{ex.p}; the general solution is obtained by multiplication by an arbitrary first integral of $\D_-$:
\begin{align*}
  p_+ = (1+u_1)F(\alpha_+), \quad \alpha_+ = u_2 + \tfrac12 u_1^2,
  \quad F>0\mbox{ arbitrary}.
\end{align*}

To solve the second equation in \eqref{ex.p}, $\D_+p_- = 2p_-/(1+u_1)+F(\alpha_+)$, we fix a $(+)$-characteristic $\{\alpha_-=\mbox{const.}\}$ and parametrize it by $u_1=t$, $u_2=t-\alpha_-$ (this gives $u_1-u_2=\alpha_-$). Since the first component of $r_\pm$ equals one, parametrizing $u_1$ in this way gives $\D_\pm=\dd/\dd u_1=\dd/\dd t$. Then the second equation of \eqref{ex.p} becomes the linear differential equation
\begin{align*}
  \frac{\dd p_-}{\dd t} = \frac{2}{1+t}p_- 
  + F\bigg(t-\alpha_-+ \frac{t^2}{2}\bigg),
\end{align*}
whose integrating factor is $(1+t)^{-2}$. This gives the solution
\begin{align*}
  p_- = (1+t)^{2}\bigg(G(\alpha_-) + \int_{t_0}^t
  \frac{F(\tau-\alpha_-+\tau^2/2)}{(1+\tau)^2}\dd\tau\bigg), \quad
  G>0\mbox{ arbitrary}. 
\end{align*}
The positivity of $p_-$ is automatic when $F$, $G>0$. 

{\em Step 3: Construction of a convex symmetrizer.}
The positive functions $F$ and $G$ are at our disposal. As an example, we choose $F=1$ and $G=c-(1+t_0)^{-1}$ for some $c\ge 1$; see below for more examples. Then $p_+=1+u_1$ and $p_- = (1+u_1)^2(G - (1+u_1)^{-1} + (1+t_0)^{-1}) = c(1+u_1)^2 - (1+u_1)$. The entropy density is recovered from Lemma \ref{lem.recon} (or from \eqref{H}):
\begin{align*}
  \na^2 h_c &= \frac{1+u_1}{(1+u_1)^2}\begin{pmatrix}
  u_1^2 & u_1 \\ u_1 & 1 \end{pmatrix}
  + \frac{c(1+u_1)^2 - (1+u_1)}{(1+u_1)^2}\begin{pmatrix}
  1 & -1 \\ -1 & 1 \end{pmatrix}
  = \begin{pmatrix}
  u_1 + c-1 & 1-c \\ 1-c & c \end{pmatrix}.
\end{align*}
When $c=1$, this can be easily integrated twice, giving $h_1(u) = u_1^3/6 + u_2^2/2$. The general case is computed from
\begin{align*}
  \na^2 h_c = (c-1)\begin{pmatrix} 1 & -1 \\ -1 & 1 \end{pmatrix}
  + \na^2 h_1.
\end{align*}
Thus,
\begin{align*}%\label{hc}
  h_c(u) = \frac{u_1^3}{6} + \frac{u_2^2}{2} + \frac{c-1}{2}(u_1-u_2)^2,
  \quad u\in\dom.
\end{align*}
The entropy density becomes separated in the variables for $c=1$, while for $c>1$, it contains the mixed term $u_1u_2$. 

\subsection{Failure mechanisms}\label{sec.fail}

Globally, the strategy presented in Section \ref{sec.main} may fail, and we briefly illustrate two possible failure mechanisms by concrete examples. 

{\em Degeneracy.} 
The strategy breaks down if the eigenbasis degenerates in the sense that $\Delta(u)$ approaches zero on some set in the interior of $\dom$; then $\delta\to 0$, and the coefficients of the transport system blow up. An example is given by the matrix
\begin{align*}
  A(u) = \begin{pmatrix}
  2 & 1 \\ \dfrac{u_2-1}{u_2+1} & 2 \end{pmatrix}
  \quad\mbox{for }u\in\dom=(0,\infty)^2.
\end{align*}
It is positively stable but $\Delta(u) = 4(u_2-1)/(u_2+1)$ may change sign. By Lemma \ref{lem.Delta}, there is no convex symmetrizer on any open subset of $\{u_2<1\}$. On the line $u_2=1$, $A(u)$ is a Jordan block, and the second clause of Lemma \ref{lem.Delta} fails there too. On the hyperbolic side $\{u_2>1\}$, we find that $\delta = m_+-m_- = \sqrt{\Delta(u)}\to 0$ as $u_2\to 1^+$, and the coefficients of the transport system blow up. Therefore, no convex symmetrizer exists on $\dom$. However, convex symmetrizers $h(u)$ exist on the whole hyperbolic side $\{u_2>1\}$. Indeed, $A$ depends only on $u_2$, and the decision test of Theorem~\ref{thm.septest} gives the separable convex symmetrizer
\begin{align*}
  h(u) = \frac{u_1^2}{2} + \frac{u_2^2}{2} + 2(u_2-1)\big(\log(u_2-1)-1\big).
\end{align*}
Notice that the blow up happens on a line in the interior of $\dom$. Replacing $u_2$ by $u_2+1$ allows us to construct a convex symmetrizer on all of $\dom$. Here, $\Delta(u)$ tends to zero on a subset of $\pa\dom$, but the blow up is fatal only if it occurs at an interior point. 

{\em No global rectangle.}
Another issue occurs when the characteristic foliations fail to produce a global rectangle $R$, so that not all of $\dom$ is covered. As an example, consider the matrix
\begin{align}\label{wind}
  A(u) = \begin{pmatrix} 
  4 + 5r^2 - 10u_1^2 & -4r^2 - 10u_1u_2 \\ 
  4r^2 - 10u_1u_2 & 4 + 5r^2 - 10u_2^2 \end{pmatrix}, \quad
  r^2 = u_1^2 + u_2^2,
\end{align}
where $u\in\dom=\{1/2<|u|<1\}$ (this domain is not simply connected). The matrix $R(\theta)=(e(\theta),e^\perp(\theta))$, where $e(\theta)=(\cos\theta,\sin\theta)^T$, $e^\perp(\theta) = (-\sin\theta,\cos\theta)^T$, describes the rotation by $\theta$. In the basis $(e(\theta),e^\perp(\theta))$, the matrix $A(u)$ becomes
\begin{align*}
  \widetilde A(r) = R(\theta)^TA(re(\theta))R(\theta)
  = \begin{pmatrix} 4-5r^2 & -4r^2 \\ 4r^2 & 4+5r^2 \end{pmatrix}.
\end{align*}
It is positively stable if $r<2/\sqrt{3}$ and $\Delta=36r^4>0$ for $r>0$. We claim that there is no convex function $h\in C^2(\dom)$ such that $\sym(\na^2 hA)$ is positive definite on $\dom$. Turning to polar coordinates, we may assume that $h(u)=\phi(r)$, where $\phi''(r)\ge 0$. Otherwise, we consider its rotational average. (Indeed, if $h$ has the required properties, then so has $\widetilde h(u) = (2\pi)^{-1}\int_0^{2\pi}h(R(\phi)u)\dd\phi$: Both conditions are convex in $\na^2h$ and invariant under rotations, because of $A(R(\phi)u)=R(\phi)A(u)R(\phi)^T$, and the strict positivity survives the averaging.) A computation shows that
\begin{align*}
  \na^2 he(\theta) &= \big(\phi''(r)e(\theta)\otimes e(\theta) 
  + r^{-1}\phi'(r)e^\perp(\theta)\otimes e^\perp(\theta)\big)e(\theta)
  = \phi''(r)e(\theta), \\
  A(u)e(\theta) &= (4-5r^2)e(\theta) + 4r^2 e^\perp(\theta).
\end{align*}
If $\sym(\na^2 hA)$ were positive definite, we obtain from
\begin{align*}
  0 < e(\theta)^T\sym(\na^2 hA)e(\theta) 
  = (\na^2 he(\theta))\cdot(Ae(\theta)) = (4-5r^2)\phi''(r)
\end{align*} 
a contradiction when $r>2/\sqrt{5}$. Thus, no such function $h$ can exist on any rotationally invariant domain containing a circle $\{|u|=r\}$ with $r>2/\sqrt{5}$, in particular on $\dom$. 

%(Interestingly, the function $h(u)=\frac12|u|^2$ is convex symmetrizer on $\{r<2/\sqrt{5}\}$, since $\sym(\na^2 hA) = \sym A = (4+5r^2)I - 10u\otimes u$ has the eigenvalues $4\pm 5r^2$, which are positive for $r<2/\sqrt{5}$.) 
The reason for the failure is that the integral curves of the eigendirections are logarithmic spirals with $r=r_0\exp(-\theta/2)$ (associated to the eigenvalue $4+3r^2$) and $r=r_0\exp(-2\theta)$ (associated to the other eigenvalue $4-3r^2$). Both wind infinitely often around the origin. There is therefore no pair of characteristic curves spanning a characteristic rectangle covering the annulus $\dom$. However, $h(u)=\frac12|u|^2$ is an entropy density on the smaller disk $\{|u|<2/\sqrt{5}\}$, since $\sym(\na^2 hA) = \sym A= (4+5r^2)I-10uu^T$ has the eigenvalues $4\pm 5r^2$, which are positive for $r<2/\sqrt{5}$. 

%\subsubsection*{Positivity loss}
%Furthermore, the solution $P_\pm$ to \eqref{star2} may approach zero, which obstructs the construction of a strictly convex function $h$. More precisely, the failure happens in case of a rectangle on which every choice of positive Goursat data produces a solution that loses positivity. We do not explore this possibility in detail. 

%%%%%%%%%%%%%%%%%%%%%%%%%%%%%%%%%%%%%%%%%%%%%%%%%%%%%%%%%%%%%%

\section{Sufficient conditions ($n=2$)}\label{sec.suff}

In this section, we give three sufficient criteria under which a convex symmetrizer exists: The curvatures $\kappa_\pm$ each have a constant sign; or the coefficients $K_\pm$ in Proposition \ref{prop.nf} are sufficiently small, independently of the sign of $\kappa_\pm$, or the characteristic curves are straight.

\subsection{Curvature sign conditions}\label{sec.sign-small}

We show first that a convex symmetrizer exists if the characteristic curves each have a constant sign on the characteristic rectangle.

\begin{theorem}[Sign condition]\label{thm.const}
Let \eqref{maincond} hold. Assume that $A$ is positively stable in $\dom$ and that $\dom$ is covered by one characteristic rectangle
$R=I_+\times I_-$, on which the coefficients of \eqref{star2} are bounded. If $R$ is unbounded, assume in addition that the total
couplings (integrals over the full sides of $R$)
\begin{align*}
  k_\pm:=\sup_{\alpha_\pm}
  \int_{I_\mp}|K_\pm(\alpha_+,\alpha_-)|\dd\alpha_\mp
\end{align*}
are finite, satisfy $k_+k_-<1$, and that these integrals converge locally uniformly in $\alpha_\pm$. If $\kappa_+$ and $\kappa_-$ each have a constant sign on $R$ (either sign, or identically zero), then a
convex symmetrizer exists on $\dom$. It may be built from any smooth, positive and bounded Goursat data.
\end{theorem}

\begin{proof}
Proposition \ref{prop.nf} implies that the coefficients
$K_\pm$ of \eqref{star2} each have a constant sign.
Choose the corner of the rectangle as follows.
We prescribe $P_+$ on $\{\alpha_-=b_0\}$ and integrate
upward in $\alpha_-$ if $K_+\ge 0$, and prescribe $P_+$
on $\{\alpha_-=b_1\}$ and integrate downward if $K_+\le 0$.
Independently, we prescribe $P_-$ on $\{\alpha_+=a_0\}$
and integrate toward increasing $\alpha_+$ if $K_-\ge 0$,
and on $\{\alpha_+=a_1\}$ integrating
toward decreasing $\alpha_+$ if $K_-\le 0$.
In each case, the relevant integral equation in
\eqref{inteq} takes the form
\begin{align*}
  P_\pm = P_\pm^{\rm data}
  + \int\mbox{(nonnegative kernel)}P_\mp \dd s.
\end{align*}
We run a Picard iteration from
$(P_+^0,P_-^0)=(P_+^{\rm data},P_-^{\rm data})$.
By induction, we have $P_\pm^n\ge P_\pm^{\rm data}>0$
for all $n\in\N$, since the integral terms are nonnegative.
For bounded $R$, the Picard iteration converges
to the unique solution by Lemma~\ref{lem.well},
and the inequalities pass to the limit.
Thus, $P_\pm\ge P_\pm^{\rm data}>0$,
and we conclude from Theorem~\ref{thm.first}.

If $R$ is unbounded, the corners may lie at infinity,
and the data are prescribed as limits.
The integral equations are then posed on intervals
of possibly infinite length.
Let $\mathcal K_\pm$ denote the corresponding
integral operators. By the assumed local uniform
convergence of the integrals of $|K_\pm|$,
these operators map $C_b(R)$ into itself and satisfy
$\|\mathcal K_\pm P\|_\infty\le k_\pm\|P\|_\infty$.
Inserting the equation for $P_+$ into the one for $P_-$
gives the fixed-point equation $P_-=\mathcal F(P_-)$, where $\mathcal F(P):=P_-^{\rm data}+\mathcal K_-P_+^{\rm data}
+\mathcal K_-\mathcal K_+P$.
Then,
\begin{align*}
  \|\mathcal F(P)-\mathcal F(Q)\|_\infty \le k_+k_-\|P-Q\|_\infty.
\end{align*}
Since $k_+k_-<1$, the second iterate of the
Picard map is a contraction. Thus, the Banach
fixed-point theorem guarantees convergence
to the unique solution.
Hence, there exists a unique bounded continuous
solution, and the positivity argument above gives
$P_\pm>0$. Finally, Theorem~\ref{thm.first} yields
a convex symmetrizer.
\end{proof}

The sign condition ensures that each characteristic family is either convex or concave throughout. The four possible corners of the rectangle correspond to the four sign combinations, which is why no particular sign is privileged. The curvatures $\kappa_\pm$ in Section \ref{sec.ex2} are $\kappa_+=0$ and $\kappa_-=-1$. 

Theorem \ref{thm.const} says nothing when a curvature changes sign. However, if the coefficients $K_\pm$ are small in a certain sense, we can still guarantee a convex symmetrizer. For this, we set
\begin{align*}
  k_+ := \sup_{\alpha_+}\int_{b_0}^{b_1}
  |K_+(\alpha_+,s)|\dd s, \quad
  k_- := \sup_{\alpha_-}\int_{a_0}^{a_1}
  |K_-(s,\alpha_-)|\dd s.
\end{align*}

\begin{theorem}[Smallness condition]\label{thm.small}
Let \eqref{maincond} hold. Assume that 
\begin{align*}
  k_+(1+2k_-)<1 \quad\mbox{and}\quad k_-(1+2k_+)<1.
\end{align*}
Then the solution to \eqref{star2} with data $P_+\equiv 1$ on $\{\alpha_-=b_0\}$ and $P_-\equiv 1$ on $\{\alpha_+=a_0\}$ satisfies $P_\pm>0$ on $R$. Consequently, a convex symmetrizer exists on $R$. In particular, $k_\pm< 1/2$ suffices. The rectangle $R$ may be unbounded (then the data are prescribed as limits).
\end{theorem}

\begin{proof}
Let $U_\pm=\sup_R|P_\pm|$. By the formulation as the integral equations \eqref{inteq}, $U_\pm$ is finite. (If $R$ is unbounded, the existence of a bounded solution follows from the fact that the assumptions imply $k_\pm<1$; hence, the right-hand side of \eqref{inteq} is a contraction with respect to the supremum norm on bounded continuous functions, with constant $\max\{k_+,k_-\}<1$.) Moreover, by \eqref{inteq} again,
\begin{align*}
  |P_+|\le 1 + \int|K_+||P_-|\dd s \le 1 + k_+U_-, \quad
  |P_-|\le 1 + \int|K_-||P_+|\dd s \le 1 + k_-U_+,
\end{align*}
from which we deduce that $U_+\le 1+k_+U_-$ and $U_-\le 1+k_-U_+$. The conditions $k_\pm(1+2k_\mp)<1$ imply that $k_\pm<1$ and $k_+k_-<1$. Substituting, $U_+\le 1 + k_+(1+k_-U_+)$ and hence
\begin{align*}
  U_+ \le \frac{1+k_+}{1-k_+k_-} \quad\mbox{and analogously,}\quad
  U_- \le \frac{1+k_-}{1-k_+k_-}.
\end{align*}
We conclude that
\begin{align*}
  P_+ \ge 1 - \int|K_+||P_-|\dd s \ge 1-k_+U_- 
  \ge 1- \frac{k_+(1+k_-)}{1-k_+k_-} > 0,
\end{align*}
and the last inequality is equivalent to $k_+(1+k_-)<1-k_+k_-$ or $k_+(1+2k_-)<1$. Similarly, $P_-\ge 1-k_-U_+>0$ if and only if $k_-(1+2k_+)<1$. 
\end{proof}

\begin{example}\rm 
We consider
\begin{align*}
  A(u) = \begin{pmatrix} 2 & 1 \\ 1-f(u_1)^2 & 2 + 2f(u_1) \end{pmatrix}
  \quad\mbox{for }u\in\dom=\R^2,
\end{align*}
where $f(u_1)=\eps/(1+u_1^2)$ and $\eps>0$. For $\eps=0$, the matrix $A(u)$ is positive definite, and a convex symmetrizer is given by $h(u)=u_1^2+u_2^2$. When $\eps>0$, the curvatures $\kappa_+ = \kappa_- = f'(u_1) = -2\eps u_1/(1+u_1^2)^2$ change sign at $u_1=0$. A computation shows that the coefficients are $K_\pm = \pm\frac14 f'(u_1)\exp(\mp f(u_1))$ and the numbers $k_\pm$ become $k_+ = 1-e^{-\eps}$, $k_- = e^\eps - 1$. Theorem \ref{thm.small} (applied on the unbounded rectangle $R=\R^2$) implies that there exists a convex symmetrizer if $\eps<\log(1+3^{-1/2})$.
 The solution to the transport system \eqref{star} is $p_\pm = 1\pm f(u_1)$, and the reconstruction Lemma \ref{lem.recon} gives $\na^2 h = \operatorname{diag}(\frac12(1-f(u_1)^2),\frac12)$ such that
\begin{align*}
  h(u) = \frac14\big(u_1^2+u_2^2 - \eps^2 u_1\arctan u_1\big)
  \quad\mbox{for }u\in\R^2.
\end{align*}
It turns out that this function is a convex symmetrizer even for all $\eps<1$. This shows that the condition in Theorem \ref{thm.small} is not sharp. 
\end{example}

In the extreme case of straight characteristics, i.e.\ $\kappa_\pm=0$, Proposition \ref{prop.nf} gives $K_\pm=0$. Then equations \eqref{star2} decouple, $\pa_{\alpha_-}P_+=\pa_{\alpha_+}P_-=0$. Thus, $P_\pm$ depends only on $\alpha_\pm$, and any positive $P_\pm^0$ gives a convex symmetrizer. The solution set is parametrized by two arbitrary functions of one variable. 

%%%%%%%%%%%%%%%%%%%%

\subsection{Straight characteristics}\label{sec.straight}

The preceding argument is formulated in characteristic
coordinates. We give a global construction for straight
characteristics directly on the state domain $\dom$, without assuming
a characteristic rectangle. 

\begin{theorem}[Straight characteristics]\label{thm.straight}
Let $\dom\subset\R^2$ be open and convex, and let $A$ be smooth
and positively stable with $\Delta>0$ in $\dom$. Assume that there
exists $\Psi\in C^\infty(\dom)$ with $\na^2\Psi>0$ and compact
sublevel sets in $\dom$. Furthermore, assume that both families of characteristic curves are straight lines (if \eqref{maincond} holds, this means that $\kappa_\pm=0$ in $\dom$). Then there exists smooth convex symmetrizer on $\dom$.
\end{theorem}

\begin{proof}
Since $\Delta>0$, $A$ has two distinct real eigenvalues.
As $\dom$ is convex, we can choose smooth unit eigenvector
fields $e_1,e_2$ on $\dom$ (they relate to $r_\pm$ from Section \ref{sec.main} only by normalization, i.e.\ $e_i=\pm r_i/|r_i|$).
Set $\D_i=e_i\cdot\na$. Since the characteristics are straight lines,
the unit vectors $e_i$ are constant along them, so that $\D_i e_i=0$.

{\em Step 1: Construction of the initial curves.}
Fix $i=1,2$, and let $\gamma$ be a maximal characteristic
of $e_i$, parametrized by $\gamma'=e_i(\gamma)$.
Since the characteristic is straight and $e_i$ is a unit
vector, we have $\gamma''=0$. Thus,
\begin{align}\label{second_d}
  \frac{\dd^2}{\dd t^2}\Psi(\gamma(t))
  &=\gamma'(t)^T\na^2\Psi(\gamma(t))\gamma'(t)
  +\na\Psi(\gamma(t))\cdot\gamma''(t)\notag\\
  &=e_i(\gamma(t))^T\na^2\Psi(\gamma(t))
  e_i(\gamma(t))>0.
\end{align}
By the compactness assumption, $\Psi(\gamma(t))\to\infty$
at both ends of the maximal parameter interval.
Since the second derivative is positive from \eqref{second_d}, there is a unique
$t_0$ such that
$\D_i\Psi(\gamma(t_0))=0$.
Thus, every characteristic of $e_i$ contains exactly one
point where $\D_i\Psi=0$. We denote the set of these points by
$\Sigma_i=\{u\in\dom:\D_i\Psi(u)=0\}$.
Since $\D_i e_i=0$, we have
$\na(\D_i\Psi)\cdot e_i = \D_i(\D_i\Psi) = e_i^T(\na^2\Psi)e_i
+ \na\Psi\cdot \D_ie_i = e^T_i(\na^2\Psi)e_i > 0$, and consequently,
$\na(\D_i\Psi)\neq0$.
Thus, by the implicit function theorem, $\Sigma_i$
is a smooth curve. For $u\in\Sigma_i$, the tangent space of $\Sigma_i$ at $u$ equals
\begin{align*}
  T_u\Sigma_i=\{v\in\R^2:\na(\D_i\Psi)(u)\cdot v=0\}.
\end{align*}
Since $\na(\D_i\Psi)(u)\cdot e_i(u)>0$, we obtain
$e_i(u)\notin T_u\Sigma_i$.
Thus, $e_i$ is transverse to $\Sigma_i$.
For each $u\in\dom$, let $t_i(u)$ be the unique number such that
$u-t_i(u)e_i(u)\in\Sigma_i$.
The implicit function theorem shows that $t_i(u)$,
and hence the intersection point, depends smoothly on $u$.

{\em Step 2: Derivation of the transport equations.}
Since $e_1,e_2$ form a basis, there are smooth functions $a,b,c,d$ such that $\D_2e_1=a e_1+b e_2$ and $\D_1e_2=c e_1+d e_2$. 
For positive functions $p_1,p_2$ to be determined, set $E=(e_1,e_2)$ and $H=E^{-T}\operatorname{diag}(p_1,p_2)E^{-1}$.
It holds that
\begin{equation}\label{eq.hessian.coeff}
  e_1^THe_1=p_1,\quad
  e_2^THe_2=p_2,\quad
  e_1^THe_2=0.
\end{equation}
In particular, $H$ is symmetric and positive definite. If $\lambda_i$
is the eigenvalue associated to $e_i$, then
\begin{align*}
  E^THAE = (E^THE)(E^{-1}AE) = \operatorname{diag}(p_1,p_2)
  \operatorname{diag}(\lambda_1,\lambda_2)
  = \operatorname{diag}(\lambda_1p_1,\lambda_2p_2),
\end{align*}
so $HA$ is symmetric as well.

It remains to ensure that $H$ is a Hessian.
Since $\dom$ is convex, this is equivalent to
$\pa_2H_{11}=\pa_1H_{12}$ and
$\pa_2H_{12}=\pa_1H_{22}$.
To derive the transport equations for $p_1$ and $p_2$,
we differentiate \eqref{eq.hessian.coeff}
along the characteristic directions.
The product rule gives
\begin{align*}
 \D_2p_1&=e_1^T(\D_2H)e_1+2ap_1,
 &0=\D_1(e_1^THe_2)&=e_1^T(\D_1H)e_2+cp_1,\\
 \D_1p_2&=e_2^T(\D_1H)e_2+2dp_2,
 &0=\D_2(e_1^THe_2)&=e_2^T(\D_2H)e_1+bp_2.
\end{align*}
Consequently, with $V=(\D_2H)e_1-(\D_1H)e_2$, we have
\begin{align*}
 \D_2p_1-(2a-c)p_1=e_1^TV,\qquad
 \D_1p_2-(2d-b)p_2=-e_2^TV.
\end{align*}
Since $e_1,e_2$ are linearly independent, $V=0$ holds if and only if
\begin{align}\label{straight.transport}
 \D_2p_1=(2a-c)p_1,\qquad
 \D_1p_2=(2d-b)p_2.
\end{align}
Writing $e_i=(e_i^1,e_i^2)^T$ and using
$\D_i=e_i^1\pa_1+e_i^2\pa_2$, we obtain
\begin{align*}
V=(\pa_1H)(e_2^1e_1-e_1^1e_2)
 +(\pa_2H)(e_2^2e_1-e_1^2e_2)
=\det(e_1,e_2)
\begin{pmatrix}
\pa_2H_{11}-\pa_1H_{12}\\
\pa_2H_{12}-\pa_1H_{22}
\end{pmatrix}.
\end{align*}
Thus, since $\det(e_1,e_2)\neq0$, system
\eqref{straight.transport} implies that 
$\pa_2H_{11}=\pa_1H_{12}$ and
$\pa_2H_{12}=\pa_1H_{22}$.

{\em Step 3: Construction of the convex symmetrizer.}
For $i=1,2$, let $\gamma_i$ be a characteristic
satisfying $\gamma_i'=e_i(\gamma_i)$ and
$\gamma_i(0)\in\Sigma_i$. Along $\gamma_2$, the first equation in
\eqref{straight.transport} becomes
\begin{align*}
\frac{\dd}{\dd t}p_1(\gamma_2(t))
=\D_2p_1(\gamma_2(t))=(2a-c)(\gamma_2(t))p_1(\gamma_2(t)).
\end{align*}
Similarly, the second equation becomes a differential equation along $\gamma_1$. The solutions to these equations become
\begin{align*}
 p_1(\gamma_2(t))
 &=\phi_1(\gamma_2(0))
   \exp\bigg(\int_0^t(2a-c)(\gamma_2(s))\dd s\bigg),\\
 p_2(\gamma_1(t))
 &=\phi_2(\gamma_1(0))
   \exp\bigg(\int_0^t(2d-b)(\gamma_1(s))\dd s\bigg).
\end{align*}
By Step~1, each characteristic of $e_i$ intersects
$\Sigma_i$ at a unique point, which depends smoothly on $u$. The segment joining $u$ to this point is a compact piece of the characteristic through $u$ and hence contained in $\dom$, so the integrals are finite. Thus, the formulas define smooth
positive functions $p_1,p_2$ on all of $\dom$.

By Step~2, the resulting matrix $H$ satisfies the
Hessian integrability conditions. Since $\dom$ is simply
connected, two applications of the Poincar\'e lemma yield
a smooth function $h$ with $\na^2h=H$.
Thus, $h$ is strictly convex, and $\na^2hA$ is positive
definite by Lemma~\ref{lem.free}.
For fixed initial data $\phi_1,\phi_2$, the function $h$
is unique up to an affine function.
\end{proof}

\begin{remark}\label{rem.straight-domains}
The assumption on $\Psi$ holds on common state domains.
For $\dom=(0,\infty)^2$, we may take
$\Psi(u)=u_1^2+u_2^2-\log u_1-\log u_2$.
For $\dom=\{u_1,u_2>0:u_1+u_2<1\}$, we may choose
$\Psi(u)=-\log u_1-\log u_2-\log(1-u_1-u_2)$.
For a general open convex domain $\dom\neq\R^2$, a natural candidate is $\Psi(u)=|u|^2-\log\operatorname{dist}(u,\pa\dom)$, which is strictly convex (since $\operatorname{dist}(\cdot,\pa\dom)$ is concave on convex domains) and has compact sublevel sets; it is not smooth in general, but it can be regularized while keeping these properties.

The convexity of $\dom$ is used only through its simple connectivity, namely for the existence of global smooth eigenvector fields and for the Poincar\'e lemma. In Step~1, each maximal characteristic is a connected piece of a straight line contained in $\dom$, and the argument applies without changes. Therefore, Theorem~\ref{thm.straight} also holds for simply connected open domains $\dom$.
\end{remark}

As an example, the Maxwell--Stefan equations of Section \ref{sec.ms2} have vanishing characteristic curvatures. Example \ref{sec.ex2} has one vanishing curvature, $\kappa_+=0$. Thus, the $(+)$-characteristics are straight, while its $(-)$-characteristics are parabolas.

\section{Separable entropies ($n=2$)}\label{sec.class}

If we restrict our search to separable strictly convex symmetrizers $h(u)=\sum_{i=1}^2 h_i(u_i)$, we can derive necessary and sufficient conditions for their existence, and the convex symmetrizer is unique (up to a positive factor and affine functions). We start by the observation that the symmetry condition for $\na^2hA$ becomes 
\begin{align}\label{Fsep}
  A_{12}(u)h_1''(u_1) = A_{21}(u)h_2''(u_2) \quad\mbox{on }\dom.
\end{align}

\begin{theorem}[Classification of convex symmetrizers]
\label{thm.class}
The following statements are equivalent:
\begin{enumerate}[label=\rm(\roman*)]
\item The matrix $A(u)$ is positively stable and admits a separable convex symmetrizer $h$.
\item There exist functions $c(u)$, $a_{12}(u_2)>0$, and $a_{21}(u_1)>0$ such that, with $a_{11}:=A_{11}$ and $a_{22}:=A_{22}$, we have $a_{11}+a_{22}>0$ and $a_{11}a_{22}>c^2 a_{12}a_{21}$ in $\dom$ and
\begin{align*}
  A(u) = \begin{pmatrix} a_{11}(u) & c(u)a_{12}(u_2) \\ 
  c(u)a_{21}(u_1) & a_{22}(u) \end{pmatrix}.
\end{align*}
\item It holds that $A(u)=M(u)^{-1}S(u)$ with $M(u)=\operatorname{diag}(a_{21}(u_1),a_{12}(u_2))$, where $a_{12}$, $a_{21}>0$, and $S(u)$ is symmetric positive definite. 
\end{enumerate}
In each case, $h_1''=a_{21}$, $h_2''=a_{12}$, and $\na^2 h(u)A(u)=S(u)$. 
\end{theorem}

\begin{proof}
(i) $\Rightarrow$ (ii): Set $a_{12}:=h_2''$, $a_{21}:=h_1''$, and $c:=A_{12}/a_{12}$. We infer from \eqref{Fsep} that $A_{21} = A_{12}h_1''/h_2''=ca_{21}$. Positive stability gives $\operatorname{tr}A = a_{11}+a_{22}>0$ and $\det A=a_{11}a_{22}-c^2a_{12}a_{21}>0$.

(ii) $\Rightarrow$ (iii): We set $M=\operatorname{diag}(a_{21},a_{12})$. Then
\begin{align*}
  S(u) = M(u)A(u) = \begin{pmatrix}
  a_{21}a_{11} & c a_{12}a_{21} \\ 
  c a_{12}a_{21} & a_{12}a_{22} \end{pmatrix}
\end{align*}
is symmetric and satisfies $\det S(u) = a_{12}a_{21}(a_{11}a_{22}-c^2a_{12}a_{21})>0$. Since $a_{11}a_{22}>c^2a_{12}a_{21}\ge 0$ and $a_{11}+a_{22}>0$ force $a_{11}$, $a_{22}>0$, the matrix $S(u)$ is positive definite.

(iii) $\Rightarrow$ (i): The function $h=\iint a_{21} + \iint a_{12}$ has the Hessian $M=\na^2 h$, and $\na^2hA=MA=S$ is symmetric positive definite, so $h$ is a separable convex symmetrizer. Then $M^{1/2}AM^{-1/2}=M^{-1/2}SM^{-1/2}$ is symmetric, positive definite. The matrices $M^{1/2}AM^{-1/2}$ and $A$ are similar, so they have the same eigenvalues, implying that $A$ is positively stable.
\end{proof}

A simple test is contained in the following proposition.

\begin{theorem}[Decision test]\label{thm.septest}
Let $\dom$ be a product of two open intervals and let $A_{12}A_{21}>0$ on $\dom$. The matrix $A$ is positively stable and admits a separable convex symmetrizer if and only if $\operatorname{tr}A>0$, $\det A>0$, and
\begin{align}\label{septest}
  \pa_{12}^2\log\frac{A_{12}}{A_{21}}=0 \quad\mbox{on }\dom.
\end{align}
In that case, the convex symmetrizer is computed, for any base point $(u_1^0,u_2^0)\in\dom$, from
\begin{align*} 
  h_1''(u_1) = \frac{A_{21}(u_1,u_2^0)}{A_{12}(u_1,u_2^0)}, \quad
  h_2''(u_2) = \frac{A_{12}(u_1^0,u_2)}{A_{21}(u_1^0,u_2)}
  \frac{A_{21}(u_1^0,u_2^0)}{A_{12}(u_1^0,u_2^0)},
\end{align*}
and the separable symmetrizer is unique up to a positive constant and an additive affine function. 
\end{theorem}

\begin{proof}
It follows from \eqref{Fsep} that
\begin{align*}
  \log\frac{h_2''(u_2)}{h_1''(u_1)} = \log\frac{A_{12}(u)}{A_{21}(u)}.
\end{align*}
Because of the separability of $h$, the mixed second derivative of the left-hand side vanishes, which yields \eqref{septest}. Conversely, let \eqref{septest} hold. Since $\dom$ is a product of two intervals, integrating $\pa_{12}^2\log(A_{12}/A_{21})=0$ shows that $\log(A_{12}/A_{21}) = g_2(u_2)-g_1(u_1)$ for some functions $g_1$, $g_2$. Then $h_1'':=\exp g_1$ and $h_2'':=\exp g_2$ are positive and satisfy $h_2''/h_1''=A_{12}/A_{21}$, i.e.\ \eqref{Fsep}, so $\na^2hA$ is symmetric. (The formulas in the theorem correspond to a particular normalization of $g_1$, $g_2$.) The positive definiteness of $\na^2 hA$ follows from the symmetry and Lemma~\ref{lem.free}, since $\na^2h=\operatorname{diag}(h_1'',h_2'')$ is positive definite and $A$ is positively stable. For uniqueness, let $h$ and $\widehat h$ be two separable convex symmetrizers. Then $h_2''/h_1''=A_{12}/A_{21}=\widehat{h}_2''/\widehat{h}_1''$, which gives $\widehat{h}_1''/h_1'' = \widehat{h}_2''/h_2''$. The left-hand side only depends on $u_1$, while the right-hand side depends only on $u_2$, thus the quotients are constant. For any $(u_1^0,u_2^0)\in\dom$, computing
\begin{align*} 
  h_1''(u_1) &= \frac{A_{21}(u_1,u_2^0)}{A_{12}(u_1,u_2^0)}h_2''(u_2^0), \\
  h_2''(u_2) &= \frac{A_{12}(u_1^0,u_2)}{A_{21}(u_1^0,u_2)}h_1''(u_1^0)=\frac{A_{12}(u_1^0,u_2)}{A_{21}(u_1^0,u_2)}
  \frac{A_{21}(u_1^0,u_2^0)}{A_{12}(u_1^0,u_2^0)}h_2''(u_2^0),
\end{align*}
and dropping the constant factor $h_2''(u_2^0)$, we finish the proof.
\end{proof}

Theorem \ref{thm.first} yields convex symmetrizers with two free functions of one variable, while separable symmetrizers come with a one-parameter family or not at all. This shows that the separability condition is rather strong.

\begin{example}[SKT model]\label{ex.skt2}\rm 
The two-species Shigesada--Kawasaki--Teramoto (SKT) model describes two population species that have the tendency to segregate from each other \cite{SKT79}. The model is given by equations \eqref{1.eq} with the positively stable diffusion matrix
\begin{align*}
  A(u) = \begin{pmatrix}
  a_{10} + 2a_{11}u_1 + a_{12}u_2 & a_{12}u_1 \\
  a_{21}u_2 & a_{20} + a_{21}u_1 + 2a_{22}u_2 \end{pmatrix},
\end{align*}
where $a_{ij} > 0$ and $u\in\dom=(0,\infty)^2$. We wish to apply the decision test of Theorem \ref{thm.septest}. Condition \eqref{septest} is fulfilled, and \eqref{Fsep} yields $a_{12}u_1h_1''(u_1) = a_{21}u_2h_2''(u_2)$. Hence, $u_ih_i''(u_i)$ is constant, showing that the separable convex symmetrizer is
\begin{align*}
  h_B(u) = \frac{1}{a_{12}}u_1(\log u_1-1) 
  + \frac{1}{a_{21}}u_2(\log u_2-1),
\end{align*}
and the separable entropy is unique up to a positive factor and an affine function by Theorem \ref{thm.septest}. 

Notice that uniqueness only holds in the class of separable entropies. For instance, let $a_{10}=a_{20}$ and $a_{11}=a_{22}=0$. The matrix 
%$A$ has the eigenvalue $\lambda=a_{10}+a_{21}u_1+a_{12}u_2$ with associated eigenvector $v=(a_{12},a_{21})$, and 
$A^T$ has the eigenvalue $\mu=a_{10}$ with associated eigenvector $w=(-a_{21},a_{12})$. Let $g$ be any strictly convex function and set $h(u) = h_B(u) + g(w\cdot u)$. Then
\begin{align*}
  \na^2 &h(u)A(u) = \na^2 h_B(u)A(u) + g''(w\cdot u)(w\otimes w)A(u) \\
  &= \na^2 h_B(u)A(u) + g''(w\cdot u)w\otimes(A^Tw) 
  = \na^2 h_B(u)A(u) + \mu g''(w\cdot u)w\otimes w
\end{align*}
is symmetric and positive definite. This shows that the set of convex symmetrizers is infinite dimensional. The Boltzmann entropy is distinguished by being the unique separable one, and it exists for all SKT parameters. It has been used for the existence analysis in \cite{ChJu04}. 
\end{example}

%%%%%%%%%%%%%%%%%%%%%%%%%%%%%%%%%%%%%%%%%%%%%%%%%%%%%%%%%%%%%

\section{Non-existence of an entropy ($n=2$)}
\label{sec.without}

We claim that equations \eqref{1.eq} with the diffusion matrix
\begin{align}\label{non.A}
  A(u) = \begin{pmatrix}
  1-(u_1-2)^2+(u_2-2)^2 & -(u_1+u_2-4)^2 \\
  (u_1-u_2)^2 & 1 + (u_1-2)^2 - (u_2-2)^2  \end{pmatrix}
\end{align}
for $u\in\dom=(0,\infty)^2$ do not possess an entropy. We have $\operatorname{tr}A=2$, $\det A=1$, and the eigenvalues are $\lambda_\pm=1$, so $A$ is positively stable in $\dom$. We can reformulate this matrix as
\begin{align}\label{non.Aw}
  A(u) = (1+r^2)I - 2(u-c)(u-c)^T + r^2 J, \quad\mbox{where }
  J = \begin{pmatrix} 0 & -1 \\ 1 & 0 \end{pmatrix},
\end{align}
$I$ is the unit matrix, $c=(2,2)^T$, and $r^2=|u-c|^2$. Since $\Delta=0$, the theory of Section \ref{sec.main} cannot be applied.

\begin{proposition}[Non-existence of entropy]\label{prop.nonex}
Let $A(u)$ be given as in \eqref{non.A}. Then there is no convex $h\in C^2(\dom)$ such that $\sym(\na^2 hA)$ is positive definite on $\dom$. 
\end{proposition}

\begin{proof}
Suppose that such a function $h$ exists and set $H=\na^2 h$. Let $u\in\dom$ with $u\neq c$ and let $w=u-c\neq 0$. It follows from formulation \eqref{non.Aw} that $A = (1+|w|^2)I-2ww^T + |w|^2J$ and consequently,
\begin{align}\label{non.HA}
  0 &< w^T(HA)w = (Hw)\cdot(Aw) 
  = (Hw)\cdot\big((1-|w|^2)w + |w|^2Jw\big) \\
  &= (1-|w|^2)w^THw + |w|^2(Hw)^T(Jw). \nonumber 
\end{align}
We formulate this inequality in terms of polar coordinates centered at $c$. Let $e(\theta) = (\cos\theta,\sin\theta)^T$, $e^\perp(\theta) = Je(\theta)$, and $g(r,\theta) = h(c+re(\theta))$ for $0<r<2$. Since $\pa_r(c+re(\theta))=e(\theta)$, $\pa_\theta(c+re(\theta)) = re'(\theta) = re^\perp(\theta)$, the chain rule gives
\begin{align*}
  & \pa_r g = \na h\cdot e(\theta), \quad
  \pa_{rr}^2 g = e(\theta)^T H e(\theta), \\
  & \pa_{\theta r}^2 g = re^\perp(\theta)^T H e(\theta)
  + \na h\cdot e^\perp(\theta) = re^\perp(\theta)^T H e(\theta)
  + r^{-1}\pa_\theta g.
\end{align*}
With $w=re(\theta)$ and $Jw=re^\perp(\theta)$, this becomes
\begin{align*}
  w^THw = r^2\pa_{rr}^2 g, \quad
  (Hw)^T(Jw) = r^2(He(\theta))^Te^\perp(\theta)
  = r\pa_{\theta r}^2 g - \pa_\theta g.
\end{align*}
We substitute these expressions into \eqref{non.HA}:
\begin{align}\label{non.g}
  (1-r^2)\pa_{rr}^2 g + r\pa_{\theta r}^2 g - \pa_\theta g > 0
  \quad\mbox{for all }\theta\in[0,2\pi),\ 0<r<2.
\end{align}
Tthe circle of radius $r$ around $c=(2,2)$ lies in $\dom$ if and only if $r<2$, and $r>1$ makes $1-r^2<0$). Therefore, we fix $1<r<2$. The functions $\pa_r g(r,\cdot)$ and $g(r,\cdot)$ are $C^1$ and $2\pi$-periodic, which implies that
\begin{align}\label{non.intg}
  \int_0^{2\pi}\pa_{\theta r}^2g\,\dd\theta = 0, \quad
  \int_0^{2\pi}\pa_\theta g\,\dd\theta = 0.
\end{align}
Therefore, an integration of \eqref{non.g} over $(0,2\pi)$ leads to
\begin{align}\label{non.con}
  (1-r^2)\int_0^{2\pi}\pa_{rr}^2 g\,\dd\theta > 0.
\end{align}
For each $\theta$, the segment $\{c+ te(\theta):0\le t\le r\}$ lies in $\dom$ since $r<2$. Thus, by assumption, $t\mapsto h(c+te(\theta))$ is convex, and its second derivative at $t=r$ is $\pa_{rr}^2g(r,\theta)\ge 0$. Then the left-hand side of \eqref{non.con} is nonpositive, which contradicts \eqref{non.con}.
\end{proof}

The proof can be understood as follows. The rotation term involving $J$ in \eqref{non.HA} vanishes after averaging around a circle, which is expressed by \eqref{non.intg}. Only the radial part remains, which for $r>1$ is incompatible with the convexity of $h$ along radii. The rotation term keeps $A$ positively stable (without it, we have $\det A = 1-r^4<0$ if $r>1$), and this term cannot be seen after the averaging.

\begin{remark}[A family of counter-examples]\label{rem.nonex}
The strategy of the previous proof can be generalized to matrices of the type (with $w=u-c$ for a fixed $c\in\R^2$)
\begin{align*}
  A(u) = d(r^2)I + r^{-2}(a(r^2)-d(r^2))ww^T + b(r^2)J, \quad r = |w|.
\end{align*}
We need $(a(r^2)-d(r^2))/r^2$ bounded to extend $A(u)$ continuously across $w=0$. The matrix \eqref{non.A} corresponds to $a(r^2) = 1-r^2$, $b(r^2)=r^2$, and $d(r^2) = 1+r^2$. A computation shows that
\begin{align*}
  Aw = a(r^2)w+b(r^2)Jw, \quad AJw = d(r^2)Jw - b(r^2)w.
\end{align*}
This implies that the matrix $A$ reads in the orthogonal basis $(w,Jw)$ as $\left(\begin{smallmatrix} a & -b \\ b & d \end{smallmatrix}\right)$. Consequently, $\operatorname{tr}A=a+d$, $\det A=ad+b^2$. If the coefficients are chosen such that $a+d>0$ and $ad+b^2>0$, the matrix $A$ is positively stable. Repeating the argument of the proof of Proposition \ref{prop.nonex} yields the inequality
\begin{align*}
  a(r_0^2)r_0^2\pa_{rr}^2 g(r_0,\theta) + b(r_0^2)
  \big(r_0\pa_{\theta r}^2 g-\pa_\theta g\big)(r_0,\theta)
  > 0 \quad\mbox{for every }\theta,
\end{align*}
where $r_0>0$ is fixed. Averaging eliminates the second term:
\begin{align*}
  a(r_0^2)r_0^2\int_0^{2\pi}\pa_{rr}^2 g(r_0,\theta)\,\dd\theta > 0.
\end{align*}
Let $a(r_0^2)<0$. If $\na^2h$ is positive semidefinite, we obtain $\pa_{rr}^2 g(r,\theta)\ge 0$ as above, which leads to a contradiction. Therefore, if $a+d>0$, $ad+b^2>0$ on $\dom$ and $a(r_0^2)<0$ for some $r_0>0$, then there is no convex $h$ such that $\na^2hA$ is positive definite on $\dom$ containing the circle $|u-c|=r_0$. 
\end{remark}

%%%%%%%%%%%%%%%%%%%%%%%%%%%%%%%%%%%%%%%%%%%%%%%%%%%%%%%%%%%%%%

\section{Examples from physics and biology ($n=2$)}\label{sec.exam}

We construct convex symmetrizers for several examples arising from applications.

\subsection{Maxwell--Stefan model}\label{sec.ms2}

The dynamics of gaseous mixtures can be described by the Maxwell--Stefan equations, which model the diffusive transport of the components of the mixture \cite{Max66,Ste71}. The model describes the evolution of the volume fractions $u_i\ge 0$, which sum up to one, $\sum_{i=0}^n u_i=1$, and $u_0$ is often identified as the water fraction. To simplify, we restrict ourselves to the case $n=2$. Then equations \eqref{1.eq} are considered with the diffusion matrix
\begin{align*}
  & A(u) = \frac{1}{a(u)}\begin{pmatrix}
  d_{02} + (d_{12}-d_{02})u_1 & (d_{12}-d_{01})u_1 \\
  (d_{12}-d_{02})u_2 & d_{01} + (d_{12}-d_{01})u_2 \end{pmatrix},
  \quad\mbox{where} \\
  & a(u) = d_{01}d_{02}u_0 + d_{01}d_{12}u_1 + d_{02}d_{12}u_2
\end{align*}
for $u\in\dom=\{u\in(0,1)^2:u_1+u_2<1\}$, where $d_{01}$, $d_{02}$, $d_{12}>0$ are diffusion coefficients. A computation shows that
\begin{align*}
  \operatorname{tr}A(u) &= a(u)^{-1}(d_{12}(u_1+u_2) 
  + d_{02}(1-u_1) + d_{01}(1-u_2))>0, \\ 
  \det A(u) &= a(u)^{-1} > 0,
\end{align*}
so $A(u)$ is positively stable in $\dom$.
%Writing
%\begin{align*}
%  A(u) = \frac{1}{a(u)}D + \frac{1}{a(u)}{u_1 \choose u_2}
%  \otimes{d_{12}-d_{02}\choose d_{12}-d_{01}}, \quad\mbox{where }
%  D = \operatorname{diag}(d_{02},d_{01}),
%\end{align*}
%a general property of the determinant yields
%\begin{align*}
%  \det A(u) &= \frac{1}{a(u)^2}\det D
%  \bigg(1 + {u_1 \choose u_2}^T D^{-1}
%  {d_{12}-d_{02}\choose d_{12}-d_{01}}\bigg) \\
%  &= d_{01}d_{02} + d_{01}(d_{12}-d_{02})u_1 + d_{02}(d_{12}-d_{01})u_2
%  = a(u)>0,
%\end{align*}
%and the trace becomes $\operatorname{tr}A(u) = a(u)^{-1}(\D_{12}(u_1+u_2) + \D_{02}(1-u_1) + \D_{01}(1-u_2))>0$. Thus, $A(u)$ is positively stable in $\dom$. 
The numbers $m_\pm$ are the roots of the quadratic equation
\begin{align*}
  (d_{12}-d_{01})u_1 m^2 + \big((d_{02}-d_{01}) + (d_{12}-d_{02})u_1
  - (d_{12}-d_{01})u_2\big)m - (d_{12}-d_{02})u_2 = 0,
\end{align*}
and a calculation leads to $\kappa_\pm = \pa_1m_\pm + m_\pm\pa_2 m_\pm = 0$, wherever these slopes are defined and the eigenvalues are real and distinct. Thus, the characteristic curves are straight lines. When the eigenvalues are real and distinct throughout $\dom$, Theorem \ref{thm.straight} and Remark \ref{rem.straight-domains} give a global convex symmetrizer without requiring a characteristic rectangle or a sign condition on $A_{12}$. 

We next construct explicit convex symmetrizers, including the cases of coincident eigenvalues. First, we apply the decision test of Theorem \ref{thm.septest}, showing that $h(u)=h_1(u_1)+h_2(u_2)$ constructed from 
\begin{align*}%\label{h1}
  h_1''(u_1) = \frac{(d_{12}-d_{02})u_2^0}{(d_{12}-d_{01})u_1}
  = \frac{c_1}{u_1}, 
  \quad h_2''(u_2) = \frac{c_2}{u_2}
\end{align*}
for some positive constants $c_1$, $c_2$ is a convex symmetrizer if $(d_{12}-d_{02})(d_{12}-d_{01})>0$. We infer that $h(u)=c_1u_1(\log u_1-1)+c_2u_2(\log u_2-1)$. However, this requires a sign condition on the coefficients $d_{ij}$ that has no physical justification. 

To build convex symmetrizers valid for all coefficients, we make the ansatz $h(u)=\phi(\ell(u))$, where $\phi\in C^2(\R)$ and $\ell(u)=\xi\cdot u+\xi_0$ for $\xi=(\xi_1,\xi_2)\in\R^2$, $\xi_0\in\R$ is an affine form. This ansatz is motivated by the fact that the characteristics are straight lines. The hyperbolic equation \eqref{1.hyper} (multiplied by $a(u)>0$, which does not change the equation) reads as $\phi''(\ell)Q_\xi(u)=0$, where
\begin{align*}
  & Q_\xi(u) = \lambda\ell(u) + \mu \quad\mbox{and} \\
  & \lambda = (d_{12}-d_{01})\xi_1 - (d_{12}-d_{02})\xi_2, \quad
  \mu = (d_{01}-d_{02})\xi_1\xi_2 - \lambda\xi_0.
\end{align*}
The function $Q_\xi$ vanishes at a point $u$ when the level line of $\ell$ through $u$ has a characteristic direction at $u$. Thus $\{\ell=0\}$ is a characteristic line if and only if $Q_\xi=0$ on $\{\ell=0\}$, which holds if and only if $\mu=0$. Taking $\mu=0$, the left-hand side of the hyperbolic equation \eqref{1.hyper} at $\phi(\ell)$ with $\ell\in\R$ equals $\lambda\phi''(\ell)\ell$. This is a constant, which can be canceled against constants of other terms (see below), if we choose $\phi''(\ell)=1/\ell$. These arguments motivate the choice $h(u) = \ell(u)(\log\ell(u)-1)$. We conclude that if $\mu=0$ then the left-hand side of \eqref{1.hyper} equals $\lambda$. The idea now is to combine various linear forms $\ell_i$ such that the weighted sum of the constants $\lambda_i=\lambda(\ell_i)$ vanishes. The procedure is formalized in the following result.

\begin{proposition}\label{prop.msent}
Let $\ell_i(u)=\xi^{(i)}\cdot u+\xi_0^{(i)}$ for $i=1,\ldots,N$ be affine forms satisfying $\mu=0$ or, more explicitly,
\begin{align}\label{kappa0}
  (d_{12}-d_{01})\xi_0^{(i)}\xi_1^{(i)} 
  - (d_{12}-d_{02})\xi_0^{(i)}\xi_2^{(i)}
  = (d_{01}-d_{02})\xi_1^{(i)}\xi_2^{(i)}
\end{align}
for all $\xi^{(i)}=(\xi_1^{(i)},\xi_2^{(i)})\in\R^2$, $\xi_0^{(i)}\in\R$. We assume that $\ell_i$ is strictly positive on $\dom$. Furthermore, let $c_1,\ldots,c_N\ge 0$ such that $\sum_{i=1}^N c_i\lambda(\ell_i) = 0$, with at least two vectors $\xi^{(i)}$ carrying $c_i>0$ and being linearly independent. Then
\begin{align*}
  h(u) = \sum_{i=1}^Nc_i\ell_i(u)(\log\ell_i(u)-1)  
\end{align*}
is a convex symmetrizer. In particular, $\na^2 hA$ is symmetric and positive definite.
\end{proposition} 

\begin{proof}
By our preceding arguments, the left-hand side of \eqref{1.hyper} equals $\sum_{i=1}^N c_i\lambda(\ell_i)$, and since this is supposed to vanish, the symmetry of $\na^2 hA$ follows. The Hessian $\na^2 h=\sum_{i=1}^N c_i\ell_i^{-1}\xi^{(i)}\otimes\xi^{(i)}$ is positive semidefinite, and it is positive definite as soon as two linearly independent vectors $\xi^{(i)}$ occur with $c_i>0$. Indeed, in this situation, we have $z^T\na^2 hz = \sum_{i=1}^N c_i\ell_i^{-1}(\xi^{(i)}\cdot z)^2>0$ for $z\neq 0$. Since $A$ is positively stable, Lemma \ref{lem.free} implies that $\na^2 hA$ is positive definite.
\end{proof}

There exist more convex symmetrizers. For instance, let $\xi^*=(d_{12}-d_{02},d_{12}-d_{01})^T$. Then $\lambda=0$, and the left-hand side of \eqref{1.hyper} with $h(u)=(\xi^*\cdot u)^2$ is the constant
\begin{align*}
  2(d_{12}-d_{01})(d_{12}-d_{02})(d_{01}-d_{02}).
\end{align*}
Thus, Proposition \ref{prop.msent} extends to
\begin{align*}
  h(u) = \sum_{i=1}^N c_i\ell_i(\log\ell_i-1)
  + \nu(\xi^*\cdot u)^2, 
\end{align*}
where $\nu\ge 0$ is determined from
\begin{align*}
  \sum_{i=1}^N c_i\lambda(\ell_i) 
  + 2\nu(d_{12}-d_{01})(d_{12}-d_{02})(d_{01}-d_{02}) = 0.
\end{align*}

\begin{example}[Two-parameter family of entropies]\rm
%\label{ex.twopara}
Let $N=3$ and take $\xi^{(1)}=(1,0)^T$, $\xi^{(2)}=(0,1)^T$, $\xi^{(3)}=(-1,-1)^T$ with $\xi_0^{(1)}=\xi_0^{(2)}=0$, $\xi_0^{(3)}=1$. Then Proposition \ref{prop.msent} gives the two-parameter family of convex symmetrizers
\begin{align*}
  h(u) = c_1u_1(\log u_1-1) + c_2u_2(\log u_2-1)
  + c_3(1-u_1-u_2)(\log(1-u_1-u_2)-1),
\end{align*}
where the numbers $c_i\ge 0$ need to satisfy $\sum_{i=1}^3 c_i\lambda(\ell_i)=0$, i.e.
\begin{align*}
  c_1(d_{12}-d_{01}) - c_2(d_{12}-d_{02})
  + (d_{01}-d_{02})c_3 = 0.
\end{align*}
A concrete example is given by $c_3=1+|d_{12}-d_{01}|+|d_{12}-d_{02}|$, $c_1=c_3+(d_{12}-d_{02})$, $c_2=c_3+(d_{12}-d_{01})$. 
\end{example}

\begin{example}[Boltzmann entropy]\rm
Choose $N=3$ and the vectors $\xi^{(i)}$ as in the previous example. We take $c_i=1$ for $i=1,2,3$. Then \eqref{kappa0} is fulfilled and $\sum_{i=1}^3 c_i\lambda(\ell_i) = 0$, so the function
\begin{align}\label{ms.hB}
  h_B(u) = u_1(\log u_1-1) + u_2(\log u_2-1) 
  + (1-u_1-u_2)(\log(1-u_1-u_2)-1)
\end{align}
is a convex symmetrizer. This is the Boltzmann entropy density for the Maxwell--Stefan equations with the volume fractions $u_1$, $u_2$, and $u_0:=1-u_1-u_2$. It was used for the existence analysis in \cite{JuSt13}. This entropy does not contain the coefficients $d_{ij}$ in its definition. 
\end{example}

We show that the Boltzmann entropy density is the only convex symmetrizer satisfying \eqref{1.hyper} for all $d_{ij}>0$ and being independent of $d_{ij}$.

\begin{proposition}%\label{prop.msunique}
Let $h\in C^3(\dom)$ be convex, satisfy \eqref{1.hyper} for all positive $d_{01}$, $d_{02}$, $d_{12}$ and being independent of $d_{01}$, $d_{02}$, $d_{12}$. Then there exists $C\ge 0$ such that $h=Ch_B$, up to an affine function, where $h_B$ is the Boltzmann entropy density defined in \eqref{ms.hB}. 
\end{proposition}

\begin{proof}
The hyperbolic equation \eqref{1.hyper} is written as
\begin{align*}
  & d_{01}E_{01} + d_{02}E_{02} + d_{12}E_{12} = 0, 
  \quad\mbox{where }E_{12} = -E_{01}-E_{02}, \\
  & E_{01} = -u_1\pa_{11}^2 h + (1-u_2)\pa_{12}^2 h, \quad
  E_{02} = -(1-u_1)\pa_{12}^2 h + u_2\pa_{22}^2 h. 
\end{align*}
By assumption, this vanishes for all $d_{ij}>0$, so each coefficient vanishes, $E_{01}=E_{02}=E_{12}=0$. Solving for $\pa_{11}^2 h$ and $\pa_{22}^2 h$, we find that
\begin{align}\label{ms.aux}
  u_1\pa_{11}^2 h = (1-u_2)\pa_{12}^2 h, \quad
  u_2\pa_{22}^2 h = (1-u_1)\pa_{12}^2 h.
\end{align}
We differentiate the first equation with respect to $u_2$, the second one with respect to $u_1$:
\begin{align*}
  u_1\pa_1(\pa_{12}^2 h) = (1-u_2)\pa_2(\pa_{12}^2 h) - \pa_{12}^2h,
  \quad
  u_2\pa_2(\pa_{12}^2 h) = (1-u_1)\pa_1(\pa_{12}^2 h) - \pa_{12}^2 h.
\end{align*}
The determinant of this linear system in the variables $\pa_1(\pa_{12}^2 h)$, $\pa_2(\pa_{12}^2 h)$ equals $u_1u_2-(1-u_1)(1-u_2)=-u_0\neq 0$, so there exists a unique solution, which is $\pa_1(\pa_{12}^2 h) = \pa_2(\pa_{12}^2h) = \pa_{12}^2 h/u_0$. We infer that
\begin{align*}
  \pa_1(u_0\pa_{12}^2 h) = -\pa_{12}^2 h + u_0\pa_1(\pa_{12}^2 h) = 0,
  \quad
  \pa_2(u_0\pa_{12}^2 h) = -\pa_{12}^2 h + u_0\pa_2(\pa_{12}^2 h) = 0.
\end{align*}
Hence, $u_0\pa_{12}^2 h = C$ on $\dom$ for some constant $C$ and, taking into account \eqref{ms.aux},
\begin{align*}
  \pa_{11}^2 h = \frac{u_0+u_1}{u_1}\pa_{12}^2 h
  = C\bigg(\frac{1}{u_1}+\frac{1}{u_0}\bigg), \quad
  \pa_{22}^2 h = \frac{u_0+u_2}{u_2}\pa_{12}^2 h
  = C\bigg(\frac{1}{u_2}+\frac{1}{u_0}\bigg).
\end{align*}
Thus, $\na^2 h=C\na^2 h_B$, and $C\ge0$ by convexity of $h$, finishing the proof.
\end{proof}

%%%%%%%%%%%%%%%%%%%

\subsection{Busenberg--Travis-type model}\label{sec.bt2}

Equations \eqref{1.eq} with the diffusion matrix
\begin{align*}
  A(u) = \begin{pmatrix} a_{11}u_1 & a_{12}u_1 \\ a_{21}u_2 & a_{22}u_2
  \end{pmatrix}\quad\mbox{for }u\in\dom=(0,\infty)^2
\end{align*}
were suggested by Busenberg and Travis \cite{BuTr83} in the context of epidemiological modeling and by Gurtin and Pipkin \cite{GuPi84} to describe populations that disperse to avoid crowding. While in these works the case $a_{ij}=a_i$ has been studied, we assume a more general situation with $a_{ij}>0$ and $a_{11}a_{22}>a_{12}a_{21}$.
 It follows from $\operatorname{tr}A>0$ and $\det A = (a_{11}a_{22}-a_{12}a_{21})u_1u_2>0$ that $A(u)$ is positively stable.

We derive families of convex symmetrizers from the hyperbolic equation \eqref{1.hyper} by observing that it can be expressed in a single variable $\sigma=u_2/u_1$. In particular, the known Boltzmann, Rao, and Lauren\c{c}ot--Matioc entropies \cite{JPZ22,LaMa23} follow directly from our unified approach. For this, we divide \eqref{1.hyper} by $u_1$:
\begin{align}\label{bt.hyper}
  E(h) := a_{12}\pa_{11}^2 h + (a_{22}\sigma-a_{11})\pa_{12}^2 h
  - a_{21}\sigma\pa_{22}^2 h = 0, \quad \sigma = \frac{u_2}{u_1}.
\end{align}
This motivates the calculation along the rays $\sigma$. We are looking for convex symmetrizers that are homogeneous of degree $k\in\N$, i.e.\ $h(\lambda u)=\lambda^k h(u)$ for $\lambda>0$. Taking $\lambda=1/u_1$, this leads to
\begin{align*}
  h(u) = u_1^k h\bigg(\frac{u}{u_1}\bigg)
  = u_1^k h\bigg(1,\frac{u_2}{u_1}\bigg) = u_1^k G(\sigma),
\end{align*}
where $G(\sigma):=h(1,\sigma)$. A computation shows that 
\begin{align*}
 & \pa_{11}^2 h = u_1^{k-2}\big(k(k-1)G-2(k-1)\sigma G'
  +\sigma^2G''\big), \\
 & \pa_{12}^2 h = u_1^{k-2}\big((k-1)G'-\sigma G''\big), \quad
 \pa_{22}^2 h = u_1^{k-2}G''.
\end{align*}
Hence, the ansatz $h(u)=u_1^k G(\sigma)$ in \eqref{bt.hyper} gives
\begin{align}
  & E(h) = u_1^{k-2}\big\{p(\sigma)G'' + (k-1)q(\sigma)G' 
  + a_{12} k(k-1)G\big\} = 0, \quad
  \mbox{where} \label{Eh} \\
  & p(\sigma) = \sigma\big((a_{12}-a_{22})\sigma
  + (a_{11}-a_{21})\big), \quad
  q(\sigma) = \big((a_{22}-2a_{12})\sigma - a_{11}\big). \nonumber 
\end{align}
Moreover,
\begin{align*}
  \pa_{11}^2 h &= u_1^{k-2}kG^{1-2/k}
  \big\{\sigma^2G^{1/k}(G^{1/k})''+(k-1)\big(\sigma (G^{1/k})'-
  G^{1/k}\big)^2\big\}, \\
  \det\na^2 h(u) &= \det\na^2(u_1^kG)
  = u_1^{2k-4}(k-1)\big(kGG'' - (k-1)(G')^2\big) \\
  &= u_1^{2k-4}(k-1)k^2G^{2-1/k}(G^{1/k})''.
\end{align*}
Thus, for $k>1$ and $G>0$, $h$ is convex on $\dom$ if and only if $G^{1/k}$ is convex on $(0,\infty)$, and in this situation $h$ is a convex symmetrizer (because of Lemma \ref{lem.free}).

Now, we distinguish the cases $k=1$ and $k>1$. Let first $k=1$. Then $E(h)=u_1^{-1}p(\sigma)G''=0$ yields $G(\sigma)=c_1+c_2\sigma$ and the affine function $h(u)=c_1u_1+c_2u_2$ with a vanishing Hessian. Therefore, we construct an entropy density in the limit $k\to 1^+$. For this, we expand formally $G_\eps = G_1 + \eps G_0 + O(\eps^2)$ for $\eps\to 0$ and define 
\begin{align*}
  h(u) = \lim_{\eps\to 0}\frac{1}{\eps}(u_1^{1+\eps}G_\eps - u_1G_1).
\end{align*}
It follows from $u_1^{1+\eps}=u_1(1+\eps\log u_1+O(\eps^2))$ that
\begin{align}\label{G0G1}
  h(u) = u_1\big(G_0(\sigma) + \log u_1 G_1(\sigma)\big).
\end{align}
Inserting $k=1+\eps$, so that $k-1=\eps$ and $k(k-1)=\eps+\eps^2$, as well as $G_\eps=G_1+\eps G_0+O(\eps^2)$ in equation \eqref{Eh} for $E(h)$, we obtain
\begin{align*}
  0 &= u_1^{1-\eps}E(h) 
  = u_1^{1-\eps}E(u_1^{1+\eps}G_\eps) 
  = p(\sigma)G_\eps'' + \eps q(\sigma)G_\eps'
  + a_{12}(\eps+\eps^2)G_\eps \\
  &= p(\sigma)G_1''
  + \eps\big(p(\sigma)G_0'' + q(\sigma)G_1' + a_{12}G_1\big)
  + O(\eps^2).
\end{align*}
Each coefficient must vanish: $pG_1''=0$ and $pG_0''=-(qG_1'+a_{12}G_1)$. The solution to the former equation is $G_1(\sigma)=c_1+c_2\sigma$. Inserting this expression into the latter equation, we arrive at
\begin{align}\label{G02}
  G_0'' = -\frac{c_2(a_{22}-a_{12})\sigma + (a_{12}c_1-a_{11}c_2)}{
  \sigma((a_{12}-a_{22})\sigma + (a_{11}-a_{21}))}.
\end{align}
This simplifies if we choose $c_1=a_{21}$ and $c_2=a_{12}$, leading to $G_0''=a_{12}/\sigma$ and after integrating twice, $G_0(\sigma)=a_{12}\sigma\log\sigma$, up to an affine function. Substituting $G_1$ and $G_0$ into \eqref{G0G1} gives the Boltzmann entropy density
\begin{align}\label{hB}
  h(u) = u_1(a_{12}\sigma\log\sigma + \log u_1(a_{21}+a_{12}\sigma))
  = a_{21}u_1\log u_1 + a_{12}u_2\log u_2.
\end{align}

Next, let $k=m\in\N$ with $m>1$ and choose $G(\sigma)=\sum_{j=0}^m \alpha_j\sigma^{m-j}$ with $\alpha_j>0$. This yields the family of Lauren\c{c}ot--Matioc entropies
\begin{align}\label{hLM}
  h(u) = \sum_{j=0}^m\alpha_j u_1^{j}u_2^{m-j}.
\end{align}
The coefficients are computed by inserting this expression into \eqref{Eh} and collecting equal monomials:
\begin{align*}
  i\big(a_{12}(i-1) + a_{22}(m-i)\big)\alpha_i
  = (m-i+1)\big(a_{11}(i-1) + a_{21}(m-i)\big)\alpha_{i-1}
\end{align*}
for $i=1,\ldots,m$. Setting $\alpha_0=1$, we can compute the coefficients $\alpha_i$ iteratively:
\begin{align*}
  \alpha_i = {m \choose i}\prod_{j=0}^{i-1}
  \frac{a_{11}j + a_{21}(m-j-1)}{a_{12}j + a_{22}(m-j-1)}>0
  \quad\mbox{for }i=1,\ldots,m.
\end{align*}
These coefficients are obtained in \cite{LaMa23} in the same way by assuming that $\na^2 hA$ is symmetric. It follows from \cite[Appendix A]{LaMa23} that the sum $h(u)$ is strictly convex and consequently, by Lemma \ref{lem.free}, $\na^2hA$ is positive definite. For $m=2$, 
\begin{align*}
  h(u) = \frac{a_{21}}{a_{22}}\bigg(\frac{a_{11}}{a_{12}}u_1^2 
  + 2u_1u_2 + \frac{a_{22}}{a_{21}}u_2^2\bigg)
\end{align*}
is the Rao entropy density.

We can construct more entropies from the approach used for $k\to 1^+$. Indeed, we have chosen particular values for $c_1$ and $c_2$ to determine $G_0$; see \eqref{G02}. Let $\ell(\sigma) := (a_{11}-a_{21})+(a_{12}-a_{22})\sigma$. Then we can formulate \eqref{G02} as
\begin{align*}
  G_0'' = \frac{a_{11}c_2-a_{12}c_1}{a_{11}-a_{21}}\frac{1}{\sigma}
  + \frac{a_{12}c_1-a_{21}c_2}{a_{11}-a_{21}}
  \frac{a_{12}-a_{22}}{\ell(\sigma)}.
\end{align*}
To simplify the arguments, let $a_{11}>a_{21}$ and $a_{12}>a_{22}$. Then $\ell(\sigma)>0$. We integrate twice to find that
\begin{align}\label{bt.G0}
  G_0 = \frac{a_{11}c_2-a_{12}c_1}{a_{11}-a_{21}}\sigma(\log\sigma-1)
  + \frac{a_{12}c_1-a_{21}c_2}{(a_{11}-a_{21})(a_{12}-a_{22})}
  \ell(\sigma)(\log\ell(\sigma)-1).
\end{align}
(If $\ell(\sigma)<0$, we may set $\ell^*(\sigma):=-\ell(\sigma)$ and obtain $-\ell^*(\sigma)(\log\ell^*(\sigma)-1)$). To obtain the entropy density $h(u)=u_1(G_0+\log u_1G_1)$, we multiply \eqref{bt.G0} by $u_1$ and use $u_1\sigma=u_2$, $u_1\ell(\sigma)=\ell\cdot u$, where $\ell=(a_{11}-a_{21},a_{12}-a_{22})^T$, such that
\begin{align*}
  u_1\sigma\log\sigma = u_2\log u_2 - u_2\log u_1, \quad
  u_1\ell(\sigma)\log\ell(\sigma) = (\ell\cdot u)\log(\ell\cdot u)
  - (\ell\cdot u)\log u_1.
\end{align*}
Combining this with the term $u_1\log u_1G_1=u_1\log u_1(c_1+c_2\sigma)$, it turns out that the coefficient of $u_2\log u_1$ cancels, and it follows that, up to a linear function,
\begin{align}\label{hBgen}
  h(u) &= u_1(G_0+\log u_1G_1) \\
  &= \beta_0(\ell\cdot u)(\log(\ell\cdot u)-1)
  + \beta_1 u_1(\log u_1 - 1) + \beta_2 u_2(\log u_2-1), \nonumber 
\end{align} 
where
\begin{align*}
  \beta_0 = \frac{a_{12}c_1-a_{21}c_2}{(a_{11}-a_{21})(a_{12}-a_{22})},
  \quad \beta_1 = \frac{a_{21}c_2-a_{22}c_1}{a_{12}-a_{22}}, \quad
  \beta_2 = \frac{a_{11}c_2-a_{12}c_1}{a_{11}-a_{21}}.
\end{align*}
We recover the Boltzmann entropy density when $c_1=a_{21}$, $c_2=a_{12}$, since $\beta_0=0$ in that case. By construction, $\na^2 hA$ is symmetric. The function $h$ is convex if $\beta_i\ge 0$ for $i=0,1,2$ (and strictly convex if, in addition, at least two of the $\beta_i$ are positive). Then Lemma \ref{lem.free} implies that $\na^2 hA$ is positive definite. Hence, the family of functions \eqref{hBgen} is a convex symmetrizer if $c_1$, $c_2\in\R$ are such that $\beta_i\ge 0$ for $i=0,1,2$. We summarize:

\begin{proposition}%\label{prop.btent}
The Busenberg--Travis model admits the following convex symmetrizers: the Boltzmann entropy density \eqref{hB}, the Lauren\c{c}ot--Matioc entropy densities \eqref{hLM}, and the logarithmic entropy densities \eqref{hBgen} if $a_{11}>a_{21}$, $a_{12}>a_{22}$, and $c_1$, $c_2$ are chosen such that $\beta_i\ge 0$ for $i=0,1,2$. A similar expression can be derived if $a_{11}<a_{21}$, $a_{12}<a_{22}$.
\end{proposition}

%%%%%%%%%%%%%%%%%

\subsection{Granular segregation model}\label{sec.sg}

Granular materials of two different components with densities $\rho_1$, $\rho_2\in[0,1]$ segregate when a partially filled cylinder is rotated about its longitudinal axis. Denoting by $u_1=\rho_1-\rho_2\in[-1,1]$ the relative density and by $u_2$ the dynamic angle between the (assumed) flat free surface of the mixture and the horizontal axis of the cylinder, the dynamics can be modeled by \eqref{1.eq} with the diffusion matrix 
\begin{align*}
  A(u) = \begin{pmatrix} \nu & u_1^2-1 \\ \gamma & 1 \end{pmatrix}
  \quad\mbox{for }u\in\dom,
\end{align*}
where $\nu$, $\gamma>0$ and $\dom=(-1,1)\times(0,2\pi)$; see \cite{ArTs92}. This matrix is positively stable. The strategy of Section~\ref{sec.main} relies on $\Delta = (\nu-1)^2-4\gamma(1-u_1^2) > 0$ in $\dom$,which in turn entails $(\nu-1)^2>4\gamma$. However, we wish to construct entropies without any restriction on $\nu$ and $\gamma$, allowing for $\Delta<0$, which is not covered by the general strategy of Section \ref{sec.main}.

To this end, we need to give up the symmetry requirement on $\na^2hA$. We are looking for separable entropies $h(u)=h_1(u_1)+h_2(u_2)$. Setting
\begin{align*}
  g_1(u_1) = (1-u_1^2)h_1''(u_1),\quad g_2(u_2) = \gamma h_2''(u_2),
\end{align*}
the matrix
\begin{align*}
  \sym(\na^2 hA) = \begin{pmatrix}
  \dfrac{\nu g_1(u_1)}{1-u_1^2} & \dfrac12(g_2(u_2)-g_1(u_1)) \\ 
  \dfrac12(g_2(u_2)-g_1(u_1)) & \dfrac{g_2(u_2)}{\gamma} \end{pmatrix}
\end{align*}
is positive definite if and only if
\begin{align*}
  4\nu g_1(u_1)g_2(u_2) > \gamma(1-u_1^2)(g_2(u_2)-g_1(u_1))^2
  \quad\mbox{for }u\in\dom.
\end{align*}
This leads to families of entropies. For instance, taking $g_1=1$ and $g_2=\theta>0$, we need to ensure that $4\nu\theta > \gamma(\theta-1)^2$, and integrating $(1-u_1^2)h_1''(u_1)=1$ and $\gamma h_2''(u_2)=\theta$ shows that 
\begin{align*}
  h(u) = \frac12\big((1+u_1)(\log(1+u_1)-1) 
  + (1-u_1)(\log(1-u_1)-1)\big) + \frac{\theta}{2\gamma}u_2^2.
\end{align*}
This function is admissible if $4\nu\theta>\gamma(\theta-1)^2$ holds, which is the case if and only if
\begin{align*}
  \theta_-<\theta<\theta_+, \quad\mbox{where }
  \theta_\pm = \frac{1}{\gamma}\big(\gamma+2\nu
  \pm 2\sqrt{\nu(\nu+\gamma)}\big).
\end{align*}
The choice $\theta=1$ gives exactly the entropy proposed in \cite{GJV03}. 

We finish with two observations. First, the dependence on $u_2$ may be nonpolynomial. Indeed, taking $g_1=1$ and $g_2(u_2)=1+a\sin u_2$ for some $a\in\R$, we find that
\begin{align*}
  h(u) = \frac12\big((1+u_1)(\log(1+u_1)-1) 
  + (1-u_1)(\log(1-u_1)-1)\big) + \frac{u_2^2}{2\gamma}
  - \frac{a}{\gamma}\sin u_2,
\end{align*}
which is strictly convex if $|a|<1$ and admissible if $1\pm a\in(\theta_-,\theta_+)$. 

Second, the product $\na^2 hA$ can be formulated as
\begin{align*}
  \na^2 h(u)A(u) = \begin{pmatrix}
  \nu/(1-u_1^2) & (\theta-1)/2 \\ (\theta-1)/2 & \theta/\gamma
  \end{pmatrix} + \frac{1+\theta}{2}J, \quad\mbox{where }
  J = \begin{pmatrix} 0 & -1 \\ 1 & 0 \end{pmatrix}.
\end{align*}
This shows that this function $h$ does not make $\na^2hA$ symmetric because of the rotational part involving $J$. In fact, no separable convex function does: For general $g_1$, $g_2>0$, the antisymmetric part of $\na^2hA$ equals $\frac12(g_1+g_2)J\neq0$. This part is invisible to the positive definiteness, since $z^TJz=0$ for $z\in\R^2$. 

%%%%%%%%%%%%%%%%%%%%%%%%%%%%%%%%%%%%%%%%%%%%%%%%%%%%%%%%%%%%%%

\section{The $n$-dimensional case}\label{sec.n}

The strategy to construct convex symmetrizers by solving a transport system along the eigendirections of the diffusion matrix generally breaks down for $n>2$. To explain this statement, we only outline the main arguments for brevity. Let $r_i$ be the eigenvectors of $A$. Recall that we can reconstruct a convex symmetrizer $h$ from $p_i = r_i^T(\na^2 h)r_i$ if $p_i>0$ (proved in Lemma \ref{lem.recon} for $n=2$). For $n=2$, the functions $p_i$ are the solutions to the transport system 
\begin{align*}
  \D_2p_1 = \Gamma_{11}p_1+\Gamma_{12}p_2, \quad
  \D_1p_2 = \Gamma_{21}p_1+\Gamma_{22}p_2,
\end{align*}
where $\Gamma_{ij}$ are some coefficients, and $\D_1p_1$, $\D_2p_2$ are free. Since $p_1$ is only constrained along $r_2$, it may be prescribed freely along a curve in the $r_1$-direction, and similarly for $p_2$. (We have used the notation $p_+=p_1$ and $p_-=p_2$ in Section \ref{sec.main}). This gives two free functions of one variable. 

If $n=3$, we have three functions $p_i=r_i^T(\na^2h)r_i$ for $i=1,2,3$ and nine derivatives $\D_ip_j$ for $i,j=1,2,3$ of which six equations prescribe $\D_ip_j$ for $i\neq j$. The three derivatives $\D_ip_i$ for $i=1,2,3$ would remain free. However, the symmetry of $(\pa_{ijk}^3h)$ with respect to all indices gives additional algebraic constraints, which are not present in the case $n=2$. They imply that, generically, $p_i$ can be written as the product $p_i = q(u)\mu_i(u)$, where the functions $\mu_i$ are determined from $A$, and one unknown function $q(u)$ remains (we leave out the details). It has three derivatives $\D_iq$, and substituting $p_i=q(u)\mu_i(u)$ into the transport equation for $\D_kp_i$ gives two equations for each derivative $\D_kq$. So each derivative of $q$ is prescribed twice, and a nonzero solution needs strong conditions on $A$, which usually fail. Even when they hold, every derivative of $q$ is fixed, and $q$ is determined up to one constant factor. In contrast to the case $n=2$, there are no free functions. 
%%%%%%%%%%%%%%%%%

\subsection{Comparison with the classification of \cite{JeKo12}}
\label{sec.jenssen}

Jenssen and Kogan \cite{JeKo12} have given a complete classification for the case $n=3$. The system for $(p_1,p_2,p_3)$ are the transport equations and the algebraic constraints, which consist of linear relations for $(p_1,p_2,p_3)$ of rank at most two \cite[Prop.~4.2]{JeKo12}. This result does not generalize to $n\ge 4$ \cite[Rem.~4.3]{JeKo12}. Assuming that the rank is constant in $\dom$, the local solution set is described as follows.

If the rank of the algebraic part is zero, there are no algebraic constraints, and the solutions depend on three arbitrary functions of one variable \cite[Theorem 3.1]{JeKo12}. This is the analogue of the case $n=2$. Jenssen and Kogan give a complete list of possibilities for the rank-one case \cite[Theorems 4.5--4.6]{JeKo12}. It may happen that no strictly convex symmetrizer exists, or the solutions depend on at most one arbitrary function of one variable and one constant, or on at most two constants. The rank-two case is the generic case discussed above with $p_i=q\mu_i$ with functions $\mu_i$ determined by $A$. If all $\mu_i$ are nonzero, which is necessary for a strictly convex symmetrizer, then the transport equations for $\D_k p_j$ with $j\neq k$ determine $\D_k q$ in terms of $q$. Hence, $q$ is unique up to a constant factor, and the convex symmetrizers form at most a one-parameter family, up to affine functions.

The rank-two case is not treated in \cite{JeKo12}, since it coincides with the case in which the conservative systems with the given eigenbasis have only the trivial fluxes $f(u)=\lambda u$ \cite[Prop.\ 4.1III]{JeKo12}, so it is irrelevant for conservation laws. For cross-diffusion systems, eigenbases of rank two may occur; an example is the thin-film solar cell model of Section \ref{sec.tfsc}. In this case, existence of a convex symmetrizer is the exception, and uniqueness (up to a factor and affine functions) holds whenever it exists. The nontrivial question is the existence itself, which motivates the explicit criteria below.

For $n\ge3$, we do not pursue the transport approach of Section~\ref{sec.main}. In the generic case, the algebraic constraints leave at most the one-parameter family (up to affine functions), so that the question is whether a convex symmetrizer exists at all rather than how to construct one. In the non-generic cases, the local structure of the solution set is described in \cite{JeKo12}, while the existence of positive solutions remains open in that work. Instead, we present in this section three results valid for any $n\ge 2$: a characterization of separable convex symmetrizers, a uniqueness result for convex symmetrizers, and a closed-orbit criterion excluding any entropy structure.

%%%%%%%%%%%%%%%%%

\subsection{Separable entropies}

The decision test of Theorem \ref{thm.septest} can be generalized to the $n$-dimensional case. Let $h(u)=\sum_{i=1}^n h_i(u_i)$ be a separable convex function in $\dom$. Then the matrix $\na^2hA$ is symmetric if and only if
\begin{align}\label{Fn}
  A_{ij}(u)h_i''(u_i) = A_{ji}(u)h_j''(u_j) \quad\mbox{for }u\in\dom
  \mbox{ for all }i\neq j.
\end{align}
We first generalize the classification result of Theorem \ref{thm.class}. 

\begin{theorem}[Classification of convex symmetrizers]
\label{thm.classn}
The following statements are equivalent:
\begin{itemize}
\item[\rm (1)] The matrix $A(u)$ is positively stable and admits a separable convex symmetrizer.
\item[\rm (2)] It holds that $A(u)=M(u)^{-1}S(u)$ with $M(u)=\operatorname{diag}(g_1(u_1),\ldots,g_n(u_n))$, where $g_i$ are positive functions of one variable, and $S(u)$ is symmetric positive definite. 
\end{itemize}
Then $h_i''=g_i$ and $\na^2hA=S$ is positive definite.
\end{theorem}

\begin{proof}
Let (1) hold. We set $M=\na^2 h$. The matrix $S=MA$ is symmetric by \eqref{Fn} and positive definite by Lemma \ref{lem.free}, which proves (2). Next, let (2) hold. Since $S$ is symmetric positive definite, $M^{1/2}AM^{-1/2} = M^{-1/2}SM^{-1/2}$ is symmetric positive definite too. As $A$ and $M^{1/2}AM^{-1/2}$ are similar matrices, they have the same eigenvalues. It follows that $A$ is positively stable. Finally, the matrix $M$ is the Hessian of the convex function $h=\sum_{i=1}^n\iint g_i$.
\end{proof}

\begin{example}[$n$-species Busenberg--Travis model]\label{ex.btn}\rm
We consider the Busenberg--Travis-type model in $n$ dimensions with matrix entries $A_{ij}(u)=u_ia_{ij}$ for $i,j=1,\ldots,n$; see Section \ref{sec.bt2} and \cite{JPZ22}. We assume that $(a_{ij})$ is symmetric and positive definite. Then $A=M^{-1}S$ with $S_{ij}=a_{ij}$ and $M=\operatorname{diag}(1/u_1,\ldots,1/u_n)$. The entropy density is computed as
\begin{align*}
  h_B(u) = \sum_{i=1}^n\int_1^{u_i}\int_1^{z}\frac{\dd y}{y}\dd z
  = \sum_{i=1}^n\big(u_i(\log u_i-1)+1\big),
\end{align*}
which corresponds to the Boltzmann entropy density. We observe that a second entropy density is given by $h_R(u)=\frac12\sum_{i,j=1}^n a_{ij}u_iu_j$ (Rao entropy), similarly as in Section \ref{sec.bt2}. Thus, the family of entropy densities $c_1h_B+c_2h_R$ depends on two parameters $c_1,c_2\ge 0$; also see \cite[Sec.\ 4]{JeKo12}.
\end{example}

The decision test (Theorem \ref{thm.septest}) has a genuinely new ingredient for $n>2$. For the following result, we set $Q_{ij}=A_{ij}/A_{ji}$ whenever $A_{ji}\neq 0$ for $i\neq j$. 

\begin{theorem}[Decision test]\label{thm.septestn}
Let $\dom$ be a product of open intervals and let $A_{ij}A_{ji}>0$ on $\dom$ for all $i\neq j$. Then there exist positive functions $g_1,\ldots,g_n$ such that $A_{ij}(u)g_i(u_i) = A_{ji}(u)g_j(u_j)$ for $u\in\dom$ and $i\neq j$ if and only if
\begin{align}\label{cocycle}
  \pa_k\log Q_{ij} = 0, \quad \pa_{ij}^2\log Q_{ij} = 0, \quad
  Q_{ij}Q_{jk}Q_{ki} = 1,
\end{align}
where the second condition is required for all $i\neq j$, and the first and third conditions for all pairwise distinct indices $i,j,k$.
The separable convex symmetrizer $h(u)=\sum_{i=1}^n h_i(u_i)$ is constructed from $h_i''=g_i>0$. If, in addition, $A$ is positively stable, then $\na^2hA$ is positive definite by Lemma~\ref{lem.free}.
\end{theorem}

The first two conditions of \eqref{cocycle} mean that each $Q_{ij}$ depends only on $(u_i,u_j)$ and factorizes, while the third condition is identified as the Kolmogorov cycle condition. It says that on the coupling graph with edges $(i,j)$, the net flow of transition rates $Q_{ij}$ along any three-state loop or cycle vanishes. When $n=2$, the third condition in \eqref{cocycle} is vacuous, which is why it does not appear in Theorem \ref{thm.septest}. 

\begin{proof}
We first show the necessity. The function $\log Q_{ij} = \log(A_{ij}/A_{ji}) = \log g_j(u_j) - \log g_i(u_i)$ depends only on $(u_i,u_j)$ and satisfies the second condition in \eqref{cocycle}. The third condition follows from
\begin{align*}
  Q_{ij}Q_{jk}Q_{ki} = \frac{A_{ij}}{A_{ji}}\frac{A_{jk}}{A_{kj}}
  \frac{A_{ki}}{A_{ik}} = \frac{g_j}{g_i}\frac{g_k}{g_j}
  \frac{g_i}{g_k} = 1.
\end{align*}
For the sufficiency, the first two conditions in \eqref{cocycle} imply that $\log Q_{ij} = \sigma_{ji}(u_j)-\sigma_{ij}(u_i)$ for some functions $\sigma_{ij}$ of a single variable. We insert this expression into the third condition of \eqref{cocycle}:
\begin{align}\label{sep.Q}
  0 = \log Q_{ij} + \log Q_{jk} + \log Q_{ki}
  = (\sigma_{ik}-\sigma_{ij})(u_i) + (\sigma_{ji}-\sigma_{jk})(u_j)
  + (\sigma_{kj}-\sigma_{ki})(u_k).
\end{align}
A sum of functions of independent variables vanishes only if each function is constant. Hence, for each fixed $i\in\{1,\ldots,n\}$, $\sigma_{ij}$ is independent of $j$ up to an additive constant, $\sigma_{ij}(u_i) = \sigma_i(u_i)+c_{ij}$ for some $c_{ij}\in\R$. Then $\log Q_{ij} = \sigma_j(u_j)-\sigma_i(u_i)+d_{ij}$, where $d_{ij}=c_{ji}-c_{ij}$ are antisymmetric constants. We infer from \eqref{sep.Q} that $d_{ij}+d_{jk}+d_{ki}=0$. We fix a base index, say one, and set $e_i=d_{1i}$. Then, setting $k=1$ and using the antisymmetry of $d_{j1}$, we have $d_{ij} = -d_{j1}-d_{1i} = d_{1j}-d_{1i} = e_j-e_i$ and consequently, $\log Q_{ij} = (\sigma_j+e_j) - (\sigma_i+e_i)$. We define $g_i(u_i):=\exp(\sigma_i(u_i)+e_i)>0$, which solves \eqref{Fn}, and $h(u)=\sum_{i=1}^n\iint g_i(u_i)$ is a convex symmetrizer.
\end{proof}

\begin{example}[$n$-species SKT model]\label{ex.sktn}\rm
The entries of the diffusion matrix associated to the $n$-species SKT model (see Example \ref{ex.skt2}) read as
\begin{align*}
  A_{ij}(u) = \delta_{ij}\bigg(a_{i0} + \sum_{k=1}^n a_{ik}u_k\bigg)
  + a_{ij}u_i, \quad i,j=1,\ldots,n,
\end{align*}
where $a_{ii}\ge 0$ and $a_{ij}>0$ for $i\neq j$. The quotients $Q_{ij}=A_{ij}/A_{ji}=a_{ij}u_i/(a_{ji}u_j)$ satisfy the first two conditions in \eqref{cocycle}, while the third one becomes
\begin{align}\label{kolmogorov}
  a_{ij}a_{jk}a_{ki} = a_{ji}a_{kj}a_{ik}
  \quad\mbox{for all pairwise distinct }i,j,k.
\end{align}
If this condition is satisfied, Theorem \ref{thm.septestn} implies the existence of a convex symmetrizer. To construct such a function, we go back to condition \eqref{Fn}, written as
\begin{align*}
  a_{ij} h_i''(u_i)u_i = a_{ji}h_j''(u_j)u_j \quad\mbox{for }i\neq j.
\end{align*}
The left-hand side depends only on $u_i$ and the right-hand side only on $u_j$; thus, both sides equal a constant $k_{ij}>0$. We have $k_{ij}=k_{ji}$ by construction. In particular, $u_ih_i''(u_i) = k_{ij}/a_{ij}$ is a constant, and it does not depend on $j$; we call it $\pi_i>0$. Hence, $h_i''(u_i) = \pi_i/u_i$. Consequently, the convex symmetrizer reads, up to affine terms, as
\begin{align*}
  h(u) = \sum_{i=1}^n h_i(u_i) = \sum_{i=1}^n\pi_i u_i(\log u_i-1).
\end{align*}
The symmetry $k_{ij}=k_{ji}$ implies the detailed-balance condition
\begin{align}\label{dbc}
  \pi_i a_{ij} = \pi_j a_{ji} \quad\mbox{for }i\neq j.
\end{align}
Conversely, given weights $\pi_i>0$ satisfying \eqref{dbc}, the Boltzmann entropy fulfills \eqref{Fn}. We have shown that the SKT model admits a strictly convex separable symmetrizer if and only if there exist $\pi_i>0$ with \eqref{dbc}, or if and only if the Kolmogorov criterion \eqref{kolmogorov} holds, in which case we may choose $\pi_1=1$ and $\pi_j=a_{1j}/a_{j1}$ for $j\neq 1$. 

The detailed-balance condition \eqref{dbc} was used in \cite{CDJ18,LeMo17} for the global existence analysis. The entropy structure of the cross-diffusion system was linked with that one of a reversible microscopic many-particle Markov process in \cite{DDD19}, and a link between local and approximating nonlocal systems was provided in \cite{DiMo24}. 
\end{example}

%%%%%%%%%%%%%%%%%%%%%%

\subsection{Uniqueness of convex symmetrizers}\label{sec.unique}

We prove that if the rank of the algebraic constraints equals $n-1$ then any convex symmetrizer is unique up to a positive factor and affine functions. In fact, we prove more: Given $m$ symmetrizers with linearly independent Hessians, the rank is at most $n-m$. Before we make this statement precise, we need some preparations. Assume that $A$ has $n$ distinct real eigenvalues in $\dom$. The algebraic constraints are given by the equations $\pa_k\pa_{ij}^2h-\pa_i\pa_{kj}^2h=0$ for $i,j,k=1,\ldots,n$. We formulate these equations for general symmetric matrices. For this, let $u\in\dom$ and let
\begin{align*}
  S(u) = \{H\in\R^{n\times n}:H^T=H,\ HA(u)=A(u)^TH\}
\end{align*}
be the space of symmetric matrices $H$ for which $HA(u)$ is symmetric. Since the eigenvalues are distinct by assumption, we have $\dim S(u)=n$. Fix $u\in\dom$ and a smooth basis $(G_1,\ldots,G_n)$ of $S$ on a neighborhood $V(u)$ of $u$, which exists since $A$ has distinct real eigenvalues. We write $H_c=\sum_{m=1}^n c_mG_m$. 

A symmetric field $H\in C^1(\dom)$ is a Hessian if and only if $\pa_k H_{ij} = \pa_i H_{kj}$ for all $i,j,k=1,\ldots,n$ (this follows after applying the Poincar\'e lemma twice). It is sufficient to consider $k<i=1,\ldots,n$, $j=1,\ldots,n$, which gives $N:=n\binom n2$ equations. We insert the expansion $H_c=\sum_{m=1}^n c_mG_m$:
\begin{align}\label{u.int}
  & \pa_k H_{c,ij} - \pa_i H_{c,kj} = (D\na c+Ec)_{kij},
  \quad\mbox{where} \\
  & (D(u)\na c(u))_{kij} := \sum_{m=1}^n
  (G_{m,ij}\pa_k c_m - G_{m,kj}\pa_i c_m)(u), \nonumber \\
  & (E(u)c(u))_{kij} := \sum_{m=1}^n(\pa_k G_{m,ij}-\pa_i G_{m,kj})(u)
  c_m(u). \nonumber 
\end{align}
We say that the {\em algebraic part} is given by
\begin{align*}
  \mathcal{A}(u) = \{E(u)^Ty: y\in\R^N,\ D(u)^Ty=0\}\subset\R^n
\end{align*}
This is the set of coefficient vectors of all linear combinations of \eqref{u.int} in which the derivatives cancel. Its dimension $r=r(u)\in\{0,\ldots,n\}$ is the {\em rank}. It holds that $\mathcal{A}(u)=\{0\}$ and $r=0$ if $n=2$, $r\le 2$ if $n=3$, and $r\le n$ if $n\ge 4$. Its kernel and the space of admissible values 
\begin{align*}
  \mathcal{K}(u) = \{H_c(u):c\in\operatorname{ker}\mathcal{A}(u)\}
  \subset S(u)
\end{align*}
have both dimension $n-r(u)$. 

\begin{lemma}\label{lem.ker}
For given $v\in V(u)$ and $c_0\in\R^n$, the following statements are equivalent:
\begin{itemize}
\item[\rm (a)] $c_0\in\operatorname{ker}\mathcal{A}(v)$.
\item[\rm (b)] There exists $Z\in\R^{n\times n}$ such that $D(v)Z+E(v)c_0=0$.
\item[\rm (c)] There exists a $C^1$ curve $c$ near $v$ such that $c(v)=c_0$ and $D(v)\na c(v)+E(v)c(v)=0$.
\end{itemize}
Consequently, $H_0\in\mathcal{K}(v)$ if and only if there exists a $C^1$ field $H$ near $v$ with values in $S$ such that $H(v)=H_0$ and $\pa_kH_{ij}(v)=\pa_iH_{kj}(v)$ for all $i,j,k$. 
\end{lemma}

\begin{proof}
The equivalence (a)$\Leftrightarrow$(b) is the Fredholm alternative in finite dimensions, since $-E(v)c_0$ lies in the range of $D(v)$ if and only if it is orthogonal to $\operatorname{ker}D(v)^T$, i.e.\ $y^TE(v)c_0=0$ for all $y$ with $D(v)^Ty=0$. The second equivalence follows from $c(w)=c_0+Z(w-v)$, so that $\na c(v)=Z$ ((b)$\Rightarrow$(c)) and the choice $Z=\na c(v)$ ((c)$\Rightarrow$(b)). The characterization of $\mathcal{K}(v)$ follows from (c) and \eqref{u.int}.
\end{proof}

The following result relates the number of symmetrizers and the rank.

\begin{proposition}\label{prop.rank}
Let $u\in\dom$ and let $A$ has $n$ distinct real eigevalues in a neighborhood of $u$. Let $h_1,\ldots,h_m$ be symmetrizers near $u$. Then $\na^2h_1(u),\ldots,\na^2h_m(u)\in\mathcal{K}(u)$ and
\begin{itemize}
\item[\rm (i)] If $r=n$ on a connected open set $V$ then every symmetrizer on $V$ is an affine function.
\item[\rm (ii)] If there is a symmetrizer $h$ near $u$ with nonsingular $\na^2h(u)$, then $r(u)\le n-1$.
\item[\rm (iii)] If $\na^2h_1(u),\ldots,\na^2h_m(u)$ are linearly independent, then $r(u)\le n-m$. 
\end{itemize}
\end{proposition}

The existence of one symmetrizer forces $r\le n-1$, and every further symmetrizer with pointwise independent Hessian lowers the maximal possible rank by one. Notice that the proposition holds for any symmetrizer, convexity is not required. 

\begin{proof}
Let $h$ be a symmetrizer on a neighborhood $V'\subset V$ of $u$ and set $H=\na^2h\in C^1(V';\R^{n\times n})$. Since $H$ and $HA$ are symmetric, $H(v)\in S(v)$ for $v\in V'$. Furthermore, $H=H_c=\sum_{m=1}^nc_mG_m$ for $c_m\in C^1(V';\R^n)$ and $\pa_kH_{ij}=\pa_{kij}^3h=\pa_{ikj}^3h=\pa_i H_{kj}$ in $V'$. We infer from \eqref{u.int} that $D(v)\na c(v)+E(v)c(v)=0$ for $v\in V'$, in particular at $v=u$. It follows from (c)$\Rightarrow$(a) in Lemma \ref{lem.ker} that $c(u)\in\operatorname{ker}\mathcal{A}(u)$, i.e.\ $\na^2h(u)=H_{c(u)}(u)\in\mathcal{K}(u)$. We apply this to $h_1,\ldots,h_m$ to obtain $\na^2h_1(u),\ldots,\na^2h_m(u)\in\mathcal{K}(u)$.

(i) If $r(u)=n$ then $\operatorname{ker}\mathcal{A}(u)=\{0\}$, which implies that $\mathcal{K}(u)=\{0\}$ and $\na^2h(u)=0$ for every symmetrizer $h$ near $u$. If $r=n$ on a connected open set $V$ and $h$ is a symmetrizer on $V$, then $\na^2h\equiv 0$ on $V$ and $h$ is affine. (ii) If $\na^2h(u)$ is nonsingular, it is a nonzero element of $\mathcal{K}(u)$. Hence, $n-r(u)=\dim\mathcal{K}(u)\ge 1$. (iii) If $\na^2h_1(u),\ldots,\na^2h_m(u)$ are linearly independent elements of $\mathcal{K}(u)$, then $n-r(u)=\dim\mathcal{K}(u)\ge m$. 
\end{proof}

\begin{remark}[Relation to Jenssen and Kogan \cite{JeKo12}]\rm 
Proposition \ref{prop.rank} is an elementary consequence of the structure of the so-called $\beta$-system of Jenssen and Kogan \cite[Sec.\ 1.4]{JeKo12}. Our algebraic part coincides with theirs up to a pointwise invertible change of unknowns ($c\leftrightarrow\beta$), and the ranks are the same. For $n=3$, the work \cite{JeKo12} analyzes in Proposition 4.2 the solution set according to the rank, which is at most two. For $n\ge 4$, only some examples are given. Up to our knowledge, the counting statements of Proposition \ref{prop.rank} for general $n$ and the uniqueness result of Theorem \ref{thm.unique} below are new.
\end{remark}

The main result of this subsection is as follows.

\begin{theorem}[Uniqueness of convex symmetrizer]\label{thm.unique}
Let $n\ge 2$, $U\subset\dom$ be connected and open and assume that $A$ has $n$ distinct real eigenvalues in $U$. Assume that the rank of the algebraic constraints $\mathcal{A}(u)$ equals $r(u)=n-1$ for all $u\in U$, and let $h$ be a convex symmetrizer on $U$ such that $\na^2h(u)$ is nonsingular for all $u\in U$. Then $h$ is unique up to a positive factor and affine functions. 
\end{theorem}

\begin{proof}
We first claim that if $h_1,h_2\in C^3(U)$ such that $\na^2h_1$ is nonsingular and $\na^2h_2=q\na^2h_1$ for a function $q$, then $q$ is constant. Indeed, we have $q=\frac{1}{n}\operatorname{tr} (\na^2h_2(\na^2h_1)^{-1})\in C^1(U)$. Let $H=\na^2 h_1$. Then $\pa_k H_{ij} = \pa_i H_{kj}$ and $\pa_k(qH_{ij}) = \pa_k(\na^2 h_2)_{ij} = \pa_i(\na^2 h_2)_{kj} = \pa_i(qH_{kj})$ imply that $(\pa_kq)H_{ij} = (\pa_iq) H_{kj}$. For fixed $j$, this means that $\na q$ is parallel to the $j$th column of $H$. If $\na q\neq 0$ at some point in $U$, all columns of $H$ are parallel there. But then $H$ has rank one, a contradiction since $n\ge 2$. 

Let $h$ and $h_*$ be two symmetrizers such that $\na^2h_*(u)\neq 0$. By Proposition \ref{prop.rank}, $\na^2h(u)$ and $\na^2h_*(u)$ belong to the line $\mathcal{K}(u)$. Hence, $\na^2h=q\na^2h_*$ for some function $q$, which is constant by the previous argument.
\end{proof}

The theorem gives a sufficient criterion for the uniqueness of the (physical) convex symmetrizer that does not require the convex symmetrizer to be known in advance; the rank is computed from the matrix $A$ alone. 
%The criterion is not necessary. Indeed, we show in Lemma \ref{lem.tfsc} below that the Boltzmann entropy is unique (up to affine functions) for the thin-film solar-cell model. In Example \ref{ex.btna} below, we argue that $r(u)=n-2$, so the uniqueness criterion of Theorem \ref{thm.unique} fails. Here, the uniqueness comes from the curvature, i.e.\ from the differential part.

\begin{example}[Busenberg--Travis model]\label{ex.btna}\rm 
For the $n$-species Busenberg--Travis model of Example \ref{ex.btn}, a computation shows that $S(u)$ has the explicit basis
\begin{align*}
  G_1 = U^{-1}, \quad G_2 = \alpha, \quad
  G_m = U^{-1}(U\alpha)^{m-1}\quad\mbox{for }m=3,\ldots,n,
\end{align*}
where $U=\operatorname{diag}(u_1,\ldots,u_n)$ and $\alpha=(a_{ij})_{i,j=1}^n$. Moreover, if $r(u)=n-2$ then
\begin{align*}
  \mathcal{A}(u) = \operatorname{span}\{e_3,\ldots,e_n\}, \quad  
  \mathcal{K}(u) = \operatorname{span}\{U,\alpha\}.
\end{align*}
Indeed, it holds that $U^{-1}=\na^2 h_B$ and $\alpha=\na^2h_R$, where $h_B$ and $h_R$ are the Boltzmann and Rao entropy densities, respectively. By Proposition \ref{prop.rank}, $U,\alpha\in\mathcal{K}(u)$, and they are linearly independent since $\alpha$ is not diagonal. As $\dim\mathcal{K}(u)=n-r(u)=2$, they span $\mathcal{K}(u)$. In the coordinate $c_m$ of $H_c=\sum_{m=1}^mc_mG_m$, they correspond to $c_1$ and $c_2$, so $\operatorname{ker}\mathcal{A}(u) = \operatorname{span}\{e_1,e_2\}$, from which we conclude that $\mathcal{A}(u)=\operatorname{span}\{e_3,\ldots,e_n\}$. The inverse direction holds if the coefficients $a_{ij}$ are not too special (we do not detail this point), which means that we have $r(u)=n-2$ for generic values of $a_{ij}$. In this situation, the rank condition of Theorem \ref{thm.unique} fails. This is not surprising, since the model admits generally two families of convex symmetrizers, namely the Boltzmann entropy and the Rao entropy. The algebraic part says that the coefficients of $G_3,\ldots,G_n$ vanish, and the Hessian $\na^2h(u)$ is spanned by $G_1$ and $G_2$ only. We infer that every convex symmetrizer is of the form $h(u)=c_1h_B(u)+c_2h_R(u) = c_1\sum_{i=1}^n u_i(\log u_i-1) + \frac12 c_2\sum_{i,j=1}^n a_{ij}u_iu_j$, up to affine functions, where $c_1,c_2\ge 0$. 
\end{example}

Proposition \ref{prop.rank} and Theorem \ref{thm.unique} reveal the picture: If the rank satisfies $r(u)=n$, there are no symmetrizers at all. If $r(u)=n-1$, there is at most one entropy (if it exists, it is unique). If $r(u)=n-m$ for $m\ge 2$, at most $m$ entropies can exist. For $n=3$, the picture is shifted, since it holds $r(u)\le 2$ \cite[Rem.\ 4.3]{JeKo12}. Then the differential part (curvatures) is needed to count the entropies. Observe that all uniqueness statements concern symmetrizers. Without symmetry, an entropy only needs $\sym(\na^2hA)$ be be positive definite, which is an open condition on the $\frac12n(n+1)$-dimensional space of symmetric matrices. There is no algebraic part, and uniqueness generally fails locally (by perturbing $h(u)$).  

%%%%%%%%%%%%%%%%%%%%%%

\subsection{Example: Thin-film solar-cell model}\label{sec.tfsc}

Thin films for solar cells can be prepared by physical vapor deposition. Here, metallic components are evaporated in a high-temperature vacuum chamber, and the incoming atoms form a thin-film layer on the substrate, used as the absorption layer of the solar cell \cite{SeBa12}. The evolution of the volume fractions of the thin-film components $u_i$ is given by \eqref{1.eq} with diffusion matrix 
\begin{align}\label{tfsc}
  A_{ij}(u) = -(a_{ij}-a_{i0})u_i\quad\mbox{for }i\neq j, \quad  
  A_{ii}(u) = \sum_{j=1,\,i\neq j}^n(a_{ij}-a_{i0})u_j + a_{i0}
\end{align}
and $i,j=1,\ldots,n$, where $a_{ij}=a_{ji}$ are positive numbers, $a_{ii}:=0$, and $u\in(0,1)^n$ satisfies $\sum_{i=1}^n u_i<1$. This model was analyzed in \cite{BaEh18,HoBu22}. 

For $n\ge3$, the cycle condition in \eqref{cocycle} fails for generic coefficients, so that, in general, no separable convex symmetrizer exists. (For $n=2$, a separable convex symmetrizer exists if $(a_{12}-a_{10})(a_{12}-a_{20})>0$.) To construct a general convex symmetrizer, we introduce the variable $u_0=1-\sum_{i=1}^n u_i$ and consider the diffusion matrix in the $n+1$ variables $u_0,\ldots,u_n\in(0,1)$ satisfying $\sum_{i=0}^nu_i=1$. In fact, the above diffusion matrix comes from the cross-diffusion equations 
\begin{align*}
  \pa_t u_i = \sum_{j=0}^n\diver\big(
  a_{ij}(u_j\na u_i-u_i\na u_j)\big), \quad i=0,\ldots,n,
\end{align*}
and the diffusion matrix with entries \eqref{tfsc} is obtained by replacing $u_0$ by $1-\sum_{i=1}^n u_i$. The diffusion matrix in the $n+1$ variables becomes $\widetilde{A}(u)=(\widetilde{A}_{ij})$ with
\begin{align*}
  \widetilde{A}_{ii} = \sum_{j=0,\,j\neq i}^n a_{ij}u_j,
  \quad \widetilde{A}_{ij} = -a_{ij}u_i \quad\mbox{for }i\neq j
\end{align*}
and $i,j=0,\ldots,n$. This matrix is not positively stable (one eigenvalue vanishes since $\sum_{i=0}^n\widetilde A_{ij}=0$). Formally applying the decision test of Theorem \ref{thm.septestn} to $\widetilde A$, by the symmetry of $a_{ij}$, the values $\log Q_{ij} = \log(\widetilde A_{ij}/\widetilde A_{ji}) = \log(u_i/u_j) + \log(a_{ij}/a_{ji})$ satisfy the three conditions in \eqref{cocycle}. Since the variables $u_0,\ldots,u_n$ are constrained by $\sum_{i=0}^nu_i=1$, this argument is only heuristic, and we verify its conclusion directly below. We infer that $h_i''(u_i)=g_i(u_i)=1/u_i$ and 
\begin{align}\label{tfsc.hB}
  h(u) = \sum_{i=0}^n u_i(\log u_i-1)
  = \sum_{i=1}^n u_i(\log u_i-1) + \bigg(1-\sum_{i=1}^n u_i\bigg)
  \bigg(\log\bigg(1-\sum_{i=1}^n u_i\bigg)-1\bigg)
\end{align}
for $u=(u_1,\ldots,u_n)\in\dom=\{u_i\in(0,1):\sum_{i=1}^nu_i<1\}$ is a convex symmetrizer associated to \eqref{tfsc}. Indeed, it is proved in \cite[Sec.~2]{BaEh18} that $\na^2 hA$ is symmetric and positive definite. We claim that in some sense, the Boltzmann entropy density \eqref{tfsc.hB} is the only convex symmetrizer, up to an affine function. 

\begin{lemma}[Uniqueness of Boltzmann entropy]\label{lem.tfsc} 
Let $n\ge 2$ and $A$ be the matrix \eqref{tfsc}. Let $h\in C^2(\dom)$ be convex. If $\na^2 hA$ is symmetric for all choices of the coefficients $a_{ij}>0$ or if $\sym(\na^2 hA)$ is positive semidefinite for all choices of $a_{ij}$, then $h = Ch_B$ for some $C\ge0$, up to an additive affine function, where $h_B$ is the Boltzmann entropy density \eqref{tfsc.hB}. 
\end{lemma}

\begin{proof}
According to Theorem \ref{thm.unique}, it is sufficient to show that the rank $r(u)$ equals $n-1$ for generic coefficients $a_{ij}$. Here, we prefer a direct (and shorter) proof. The idea is to test the Hessian $\na^2 h$ against $n$ special choices of the coefficients $a_{ij}$, to conclude that this implies that $\na^2 h = C\na^2 h_B$, and finally to show that $C$ is constant. For this, fix $k\in\{1,\ldots,n\}$ and consider the matrix $A$ with all $a_{ij}=1$ except $a_{k0}=1+\theta$ (with $\theta>-1$). Then $A = I+\theta E_k$ with $E_k = u_0e_k\otimes e_k + u_k e_k\otimes\mathbf{1}$, where $e_k$ is the standard Euclidean basis vector and $\mathbf{1}$ is the vector with entries one. Let $h$ satisfy one of the two hypotheses and set $H=\na^2 h$.

If $HA$ is supposed to be symmetric, we take $\theta=1$. We deduce from the symmetry of $HA = H + HE_k$ and $H$ that $HE_k$ is symmetric too. On the other hand, if $\sym(HA) = \sym H + \theta\sym(HE_k)$ is supposed to be positive semidefinite, we divide by $\theta>0$ and pass to the limit $\theta\to\infty$ to find that $\sym(HE_k)$ is positive semidefinite. %These statements holds for all $k=1,\ldots,n$. 

Let $w_k = u_0e_k + u_k\mathbf{1}\neq 0$ and set $v=He_k$. Then $HE_k=v\otimes w_k$. We distinguish again the two cases. If $HE_k$ is symmetric, we compute $ v|w_k|^2 = (v\otimes w_k)\cdot w_k = (w_k\otimes v)\cdot w_k = w_k(v\cdot w_k)$, implying that $v=cw_k$ with $c=(v\cdot w_k)/|w_k|^2$. If $\sym(HE_k)$ is positive semidefinite and $v$, $w_k$ were linearly independent, there would be a vector $z$ with $z\cdot v=1$ and $z\cdot w_k=-1$. Then $z^T\sym(HE_k)z = z^T(v\otimes w_k)z = (z\cdot v)(z\cdot w_k)=-1<0$, contradiction. Thus, $v$ and $w_k$ are linearly dependent, and there exists $c\in\R$ such that $v=cw_k$. It follows from $z^T\sym(HE_k)z = c(z\cdot w_k)^2\ge 0$ that $c\ge 0$. Therefore, in both cases, we can write $v=c_kw_k$ for some $c_k\in\R$ (in the first case) and $c_k\ge 0$ (in the second case).

We obtain $He_k = v = c_kw_k$ or $H_{jk} = c_k(u_0\delta_{jk}+u_k)$ for $j,k=1,\ldots,n$. Since $H$ is symmetric, $c_ku_k = H_{jk} = H_{kj} = c_ju_j$ for $j\neq k$, so $\mu:=c_ju_j$ for all $j=1,\ldots,n$. We insert $c_j=\mu/u_j$:
\begin{align*}
  H_{kk} = \mu\bigg(\frac{u_0}{u_k}+1\bigg)
  = \mu u_0\bigg(\frac{1}{u_k}+\frac{1}{u_0}\bigg), \quad
  H_{jk} = \mu = (\mu u_0)\frac{1}{u_0}\quad\mbox{for }j\neq k.
\end{align*}
We infer that $H$ equals, up to the factor $\mu u_0$, the Hessian of the Boltzmann entropy density $h_B$, defined in \eqref{tfsc.hB}, $H=C\na^2 h_B$ with $C:=\mu u_0$. 

We claim that the factor $C$ is constant. The previous arguments used the symmetry of $H$ at single points. Now, we use the fact that $H$ is the Hessian of $h$. To simplify, let $h\in C^3$. (If $h\in C^2$, the arguments below hold in the sense of distributions.) Then $\pa_k H_{ij} = \pa_{ijk}^3 h = \pa_i H_{kj}$ for all $i,j,k=1,\ldots,n$. We differentiate
\begin{align*}
  \pa_k H_{ij} &= \pa_k(C\na^2 h_B)_{ij} 
  = \pa_kC(\na^2 h_B)_{ij} + C\pa_{k}\pa_{ij}^2 h_B, \\
  \pa_i H_{kj} &= \pa_i(C\na^2 h_B)_{kj} 
  = \pa_iC(\na^2 h_B)_{kj} + C\pa_i\pa_{kj}^2 h_B.
\end{align*}
Since the left-hand sides and the last term on the right-hand sides coincide, we conclude that 
\begin{align}\label{tfsc.aux}
  \pa_kC(\na^2 h_B)_{ij} = \pa_iC(\na^2 h_B)_{kj}
  \quad\mbox{for all }i,j,k=1,\ldots,n.
\end{align}
Fix $j$ and let $g:=\na C=(\pa_1C,\ldots,\pa_nC)$, $c_i=\pa_{ij}^2h_B$. Then \eqref{tfsc.aux} reads as $g_kc_i = g_ic_k$ for all $i$, $k$. As $\na^2h_B$ is positive definite, $c\neq 0$, and there exists $i$ such that $c_i\neq 0$. Then $g_k = (g_i/c_i)c_k$ for every $k$ or $g = \lambda_j c$ with $\lambda_j = g_i/c_i$ (recall that we fixed $j$). It follows that $\na C$ is a multiple of every column of $\na^2 h_B$. 

Suppose that $\na C\neq 0$ at some point. Then $\lambda_j\neq 0$ for every $j$, and every column $(\na^2 h_B)e_j=\na C/\lambda_j$ is a multiple of the vector $\na C$. It follows that $\na^2 h_B$ is of rank one, which contradicts the positive definiteness of $\na^2 h_B$. Hence, $\na C=0$ at every point, and $C$ is constant, showing the claim. Finally, $C\ge0$ since $h$ is convex.
\end{proof}

The logical structure of the proof is as follows: The special coefficient choices are test cases for one fixed function $h$; the conclusion $h = Ch_B$ is a statement about $h$, not about the coefficients. If $h$ works for all coefficients, it works in particular for the test choices, and these already force $h = Ch_B$. For a {\em fixed} choice of coefficients, the question is different. Indeed, the ansatz $h(u)=\sum_{i=0}^n c_iu_i(\log u_i-1)$ gives conditions on the coefficients $c_i$, depending on $a_{ij}$, and $c_i$ may exist or not. For instance, if $a_{i0}=:a_0$ for all $i=1,\ldots,n$ and $a_{ij}\neq a_0$ for $i,j\ge1$, $i\neq j$, one can show that 
\begin{align*}
  h(u) = c_0u_0(\log u_0-1) + c\sum_{i=1}^n u_i(\log u_i-1)
  \quad\mbox{for }c_0,c>0
\end{align*}
is a two-parameter family of convex symmetrizers. For these coefficients, $h_B$ is not the only convex symmetrizer. 

%%%%%%%%%%%%%%%%%%%

\subsection{Non-existence of an entropy}

We prove that if there exists a closed orbit for some column field from $A(u)$, then there does not exist an entropy structure.

\begin{theorem}[Closed-orbit obstruction]\label{thm.coo}
Let $w\in\R^n$, $w\neq 0$ and define $X_w:\dom\to\R^n$, $X_w(u)=A(u)w$. Suppose that there exist $T>0$ and a $C^1$ curve $\gamma:[0,T]\to\dom$ such that
\begin{align*}
  \gamma'(t) = X_w(\gamma(t))\quad\mbox{for }t\in[0,T], \quad
  \gamma(T)=\gamma(0).
\end{align*}
Then there is no function $h\in C^2(\dom)$ such that $\na^2hA$ is positive definite. In particular, there is no convex function with this property.
\end{theorem}

\begin{proof}
Suppose, by contradiction, that there exists some function $h\in C^2(\dom)$ such that $\na^2hA$ is positive definite, and let $\gamma:[0,T]\to\dom$ be a given solution curve of $\gamma'=X_w(\gamma)$ with $\gamma(T)=\gamma(0)$. We define the directional derivative function $\phi_w:\dom\to\R$ by $\phi_w(u)=\na h(u) w$. Then
%\begin{align*}
%  (\na^2h(u)w)_i = \sum_{j=1}^n\pa_{ij}^2h(u)w_j
%  = \pa_i\bigg(\sum_{j=1}^n w_j\pa_j h(u)\bigg) = \pa_i\phi_w(u),
%\end{align*}
$\na\phi_w(u)=\na^2h(u)w$. Hence, as  $\phi_w\in C^1(\dom)$,
\begin{align*}
  \frac{\dd}{\dd t}\phi_w(\gamma(t))
  &= \na\phi_w(\gamma(t))\cdot\gamma'(t)
  = \na\phi_w(\gamma(t))\cdot X_w(\gamma(t)) \\
  &= (\na^2h(\gamma(t))w)\cdot(A(\gamma(t))w) 
  = w^T\big(\na^2h(\gamma(t))A(\gamma(t))\big)w > 0.
\end{align*}
Now, on the one hand $\phi_w(\gamma(T)) - \phi_w(\gamma(0)) = 0$ by assumption, and on the other hand,
\begin{align*}
  \phi_w(\gamma(T)) - \phi_w(\gamma(0))
  = \int_0^T\frac{\dd}{\dd t}\phi_w(\gamma(t))\dd t > 0;
\end{align*}
contradiction. Hence no such function $h$ can exist.
\end{proof}

In two dimensions $n=2$, the obstruction of Theorem \ref{thm.coo} cannot happen, i.e., no map $u\mapsto A(u)w$ (for some $w\neq 0$) has a periodic orbit in $\dom$. Indeed, positive stability makes $A(u)$ invertible, so $A(u)w\neq 0$, i.e., the map has no zeros. In the plane, a periodic orbit must enclose a zero of the map (index theory), and on a simply connected $\dom$ the enclosed region lies in $\dom$; contradiction.

\begin{example}[Non-existence of entropy]\label{ex.nonex}\rm
Consider the matrix
\begin{align*}
  A(u) = \begin{pmatrix}
  1 & 0 & 1-u_2 \\ 0 & 1 & u_1-1 \\ u_2-1 & 1-u_1 & (\rho(u)-1/4)^2
  \end{pmatrix}\quad\mbox{for }u\in\dom=(0,\infty)^3,
\end{align*}
where $\rho(u)=(u_1-1)^2+(u_2-1)^2$. We compute the characteristic polynomial of $-A(u)$:
\begin{align*}
  p(\lambda) = \det(\lambda I+A(u)) = \lambda^3 + (2+a)\lambda^2
  + (1+2a+\rho)\lambda + (a+\rho),
\end{align*}
where $a=(\rho(u)-\frac14)^2$. The eigenvalues of $\sym A$ are $\lambda=1$, $\lambda=(\rho(u)-\frac14)^2$, and the latter one vanishes on the cylinder $Z=\{u\in\R^3:\rho=\frac14\}$. All coefficients of $p$ are positive (note that $a$ and $\rho$ cannot vanish simultaneously), and $(2+a)(1+2a+\rho)-(a+\rho)=2+4a+\rho+2a^2+a\rho>0$. By the Routh--Hurwitz criterion, all roots of $p$ have negative real parts, so $A$ is positively stable on $\dom$, including $Z$. The degeneracy of $\sym A$ on $Z$ is unavoidable. Indeed, if $\sym A$ is positive definite on $\dom$ then the quadratic entropy density $h(u)=|u|^2/2$ satisfies $\sym(\na^2 hA)=\sym A>0$. Thus, any non-existence result must allow $\sym A$ to fail to be positive definite somewhere. 

We define for $c>0$ the orbits $\gamma_c:[0,2\pi]\to\R^3$ by
\begin{align*}
  \gamma_c(t) = \bigg(1+\frac12\cos t,1+\frac12\sin t,c\bigg),
  \quad t\in[0,2\pi],
\end{align*}
and set $w=(0,0,1)^T$. Clearly, $\gamma_c(2\pi)=\gamma_c(0)$. The curve describes a circle of radius $\frac12$ centered at $(1,1,c)^T$ in the horizontal plane $\{u_3=c\}$. It remains to compute $X_w(u)=A(u)w=(1-u_2,u_1-1,(\rho(u)-\frac14)^2)^T = (-\frac12\sin t,\frac12\cos t,0)^T$ and
\begin{align*}
  \gamma_c'(t) = \bigg(-\frac12\sin t,\frac12\cos t,0\bigg)
  = X_w(\gamma_c(t)).
\end{align*}
Thus, Theorem \ref{thm.coo} applies, and no entropy exists.
\end{example}

%%%%%%%%%%%%%%%%%%%%%%%%%%%%%%%%%%%%%%%%%%%%%%%%%%%%%%%%%%%%%%

\end{document}